\documentclass[10pt,a4paper]{amsart}
\RequirePackage[OT1]{fontenc}
\RequirePackage{amsthm,amsmath,amssymb,mathrsfs, enumerate}
\RequirePackage{hyphenat}
\newtheorem{theorem}{Theorem}[section]
\newtheorem{corollary}[theorem]{Corollary}
\newtheorem{lemma}[theorem]{Lemma}
\newtheorem{prop}[theorem]{Proposition}
\newtheorem{definition}[theorem]{Definition}
\newtheorem{example}[theorem]{Example}
\newtheorem{probl}[theorem]{Problem}
\newtheorem{remark}[theorem]{Remark}
\newcommand{\Ii}{1\!\!1}
\usepackage{graphicx, color}
\def\RR{\mathbb{R}}
\def\NN{\mathbb{N}}
\def\r{{\bf r}}\def\r{{\bf r}}
\def\y{{\bf y}}\def\z{{\bf z}}\def\x{{\bf x}}
\def\C{{\mathcal C}}

\begin{document}
\title[The full infinite dimensional moment problem on semi-algebraic sets]{The full infinite dimensional moment problem on semi-algebraic sets of generalized functions}
\author[M. Infusino, T. Kuna, A. Rota]{M. Infusino$^{*,+}$, T. Kuna$^{+}$, A. Rota$^{+}$}
\address[Authors' address]{Department of Mathematics and Statistics, University of Reading,
Whiteknights, PO Box 220, Reading RG6~6AX, United Kingdom}
\email[Maria Infusino]{infusino.maria@gmail.com}
\email[Tobias Kuna]{t.kuna@reading.ac.uk}
\email[Aldo Rota]{a.rota@pgr.reading.ac.uk}
\thanks{$^{*}$ Supported by a Marie Curie fellowship of the Istituto Nazionale di Alta Matematica (INdAM)}
\thanks{$^{+}$ Supported by the EPSRC Research Grant EP/H022767/1}

\begin{abstract}
We consider a generic basic semi-algebraic subset~$\mathcal{S}$ of the space of generalized functions, that is a set given by (not necessarily countably many) polynomial constraints. We derive necessary and sufficient conditions for an infinite sequence of generalized functions to be \emph{realizable} on~$\mathcal{S}$, namely to be the moment sequence of a finite measure concentrated on~$\mathcal{S}$. Our approach combines the classical results about the moment problem on nuclear spaces with the techniques recently developed to treat the moment problem on basic semi-algebraic sets of~$\RR^d$. In this way, we determine realizability conditions that can be more easily verified than the well-known Haviland type conditions. Our result completely characterizes the support of the realizing measure in terms of its moments. As concrete examples of semi-algebraic sets of generalized functions, we consider the set of all Radon measures and the set of all the measures having bounded Radon-Nikodym density w.r.t.\! the Lebesgue measure.
\end{abstract}
\maketitle

\section*{Introduction}
It is often more convenient to consider characteristics of a random distribution instead of the random distribution itself and try to extract information about the distribution from these characteristics. In this paper, we are more concretely interested in distributions on functional objects like random fields, random points, random sets and random measures. The characteristics under study are polynomials of these objects like the density, the pair distance distribution, the covering function, the contact distribution function, etc.. This setting is considered in numerous areas of applications: heterogeneous materials and mesoscopic structures~\cite{To02}, stochastic geometry~\cite{Mol05}, liquid theory~\cite{HaMcDo87}, spatial statistics~\cite{Stoy00}, spatial ecology~\cite{MDL04} and neural spike trains~\cite{BKM04, JM04}, just to name a few. 

The subject of this paper is the full power moment problem on a pre-given subset $\mathcal{S}$ of $\mathscr{D}'(\RR^d)$, the space of all generalized functions on $\RR^d$. This framework choice is mathematically convenient and general enough to encompass all the aforementioned applications. More precisely, our paper addresses the question of whether certain prescribed generalized functions are in fact the moment functions of some finite measure concentrated on $\mathcal{S}$. If such a measure does exist, it will be called \emph{realizing}. The main novelty of this paper is to investigate how one can read off support properties of the realizing measure directly from positivity properties of its moment functions. 

To be more concrete, homogeneous polynomials are defined as powers of linear functionals on $\mathscr{D}'(\RR^d)$ and their linear continuous extensions. We denote by~$\mathscr{P}_{\C_c^\infty}(\mathscr{D}'(\RR^d))$ the set of all polynomials on $\mathscr{D}'(\RR^d)$ with coefficients in $\C_c^\infty(\RR^d)$, which is the set of all infinitely differentiable functions with compact support in~$\RR^d$. 

In this paper, we try to find a characterization via moments of measures concentrated on basic \textit{semi-algebraic} subsets of $\mathscr{D}'(\RR^d)$, i.e.\! sets that are given by polynomial constraints and so are of the following form
\[
\mathcal{S} = \bigcap_{i \in Y} \left\{ \left. \eta \in \mathscr{D}'(\RR^d)  \right| \ P_i(\eta) \geq 0\right\},
\]
where $Y$ is an arbitrary index set (not necessarily countable) and each $P_i$ is a polynomial in~$\mathscr{P}_{\C_c^\infty}(\mathscr{D}'(\RR^d))$. Equality constraints can be handled using $P_i$ and $-P_i$ simultaneously.
As far as we are aware, the infinite dimensional moment problem has only been treated in general on affine subsets \cite{BS71, BeKo88} and cones \cite{S74} of nuclear spaces (these results are stated in Subsection~\ref{Sec-MPNuclSpac} and Subsection \ref{ExRadonMeas}). Special situations have also been handled; see e.g. \cite{Ze83, BeKoKuLy99, KoKuOl02}. 

\subsection*{Previous results}\ \\
Characterization results via moments are built up out of five completely different types of conditions
\begin{enumerate}[I.]
\item positivity conditions on the moment sequence;
\item\label{condquasi} conditions on the asymptotic behaviour of the moments as a sequence of their degree; 
\item properties of the putative support of the realizing measure;
\item\label{condreg} regularity properties of the moments as generalized functions;
\item\label{condgro} growth properties of the moments as generalized functions. 
\end{enumerate} 
Conditions of type \ref{condreg} and \ref{condgro} are only relevant for the infinite dimensional moment problem. The general aim in moment theory is to construct a solution which is as weak as possible w.r.t.\!\! some combination of the above different types of conditions, since it seems unfeasible to get one solution which is optimal in all types simultaneously.

Let us give a review of some previous results on which our approach is based and describe the different types of conditions involved in each of them.\\
Given a sequence $m$ of putative moments, one can introduce on the set of all polynomials the so-called Riesz functional $L_m$, which associates to each polynomial its putative expectation. If a polynomial $P$ is non-negative on the prescribed support~$\mathcal{S}$, then a necessary condition for the realizability of $m$ on $\mathcal{S}$ is that $L_m(P)$ is non-negative as well. The question whether this condition alone is also sufficient for the existence of a realizing measure concentrated on $\mathcal{S}\subseteq\RR^d$ is answered by the Riesz-Haviland theorem \cite{Riesz23, Hav36}; for infinite dimensional versions of this theorem see e.g. \cite{Le75a, Le75b, Me77} for point processes and \cite{Ku09, MolLach} for the truncated case. The disadvantage of this type of positivity condition is that it may be rather difficult and also computationally expensive to identify all non-negative polynomials on~$\mathcal{S}$, especially if the latter is geometrically non-trivial.

A classical result shows that all non-negative polynomials on $\mathbb{R}$ can be written as the sum of squares of polynomials (see \cite{PolSze76}). Hence, it is already sufficient for realizability on $\mathcal{S}=\mathbb{R}$ to require that $L_m$ is non-negative on squares of polynomials, that is, $m$ is \textit{positive semidefinite}. For the moment problem on $\mathcal{S}=\RR^d$ with $d\geq 2$, the positive semidefiniteness of $m$ is no longer sufficient, as already pointed out by D.~Hilbert in the description of his 17th problem. However, the positive semidefiniteness of $m$ becomes sufficient if one additionally assumes a condition of type \ref{condquasi}, that is, a bound on a certain norm of the $n-$th putative moment $m^{(n)}$. For example, one could require that $|m^{(n)}|$ does not grow faster than $BC^n n!$ or than $BC^n \left( n \ \ln(n)\right)^n$ for some constants $B, C >0$. The weakest known growth condition of this kind is that the sequence $m$ is quasi-analytic (see\!~Appendix~\ref{Sec-App}). We will call such a sequence \textit{determining}, because this property guarantees the uniqueness of the realizing measure. The determinacy condition in the infinite dimensional case additionally involves the types~\ref{condreg} and~\ref{condgro}.
 
Beyond the results for $\mathcal{S}=\RR^d$, for a long time the moment problem was only studied for specific proper subsets~$\mathcal{S}$ of~$\RR^d$ rather than general classes of sets. However, enormous progress has recently been made for the moment problem on general basic semi-algebraic sets of~$\RR^d$. Let us mention just a few key works which were inspiring for the results presented here; for a more complete overview see~\cite{LasBook, LauMo, MarshBook}.\\
The common feature of these works is that the support properties of the realizing measure are encoded in a positivity condition stronger than the positive semidefiniteness; namely, the condition that $L_m$ is non-negative on the \textit{quadratic module} generated by the polynomials $(P_i)_{i \in Y}$ defining the basic semi-algebraic set $\mathcal{S}$. This module is the set of all polynomials given by finite sums of the form $\sum_{i} Q_i P_i$, where $Q_i$ is a sum of squares of polynomials. Semidefinite programming allows an efficient numeric treatment of such positivity conditions; see e.g.\! \cite{LasBook}. In 1982, C.~Berg and P.~H.~Maserick showed in \cite{BeMa} that for a compact basic semi-algebraic $\mathcal{S}\subset \mathbb{R}$ the positivity condition involving the quadratic module is also sufficient. Concerning the higher dimensional case, a few years later K.~Schm\"udgen proved in his seminal work~\cite{Schm91} that for a compact basic semi-algebraic $\mathcal{S}\subset\mathbb{R}^d$ a slightly stronger positivity condition, that is, $L_m$ is non-negative on the \emph{pre-ordering} generated by $(P_i)_{i \in Y}$, is sufficient. This result was soon refined by M. Putinar in~\cite{Put93} for Archimedean quadratic modules. Since then, the problem to extend their results to wider classes of~$\mathcal{S}$ has intensively been studied (see e.g.\! \cite{PowSch01, KuMarSch05, CimMarNet11}). By additionally assuming a growth condition of the type discussed above, J.~B.~Lasserre has recently showed in~\cite{Las2011} that the non-negativity of $L_m$ on the quadratic module is sufficient for realizability on a general basic semi-algebraic set $\mathcal{S}\subseteq\RR^d$.

Using the central idea of these works, we prove in this paper that also for a moment problem on an infinite dimensional basic semi-algebraic set $\mathcal{S}$, the non-negativity of $L_m$ on the associated quadratic module is sufficient for realizability under an appropriate growth condition on the sequence $m$. 

\subsection*{Outline of the contents}\ \\
Let us outline the contents and the contributions of this paper. 

In Section~\ref{Sec-Prel}, we recall some preliminaries about generalized functions, which are particularly relevant for this paper, and we pose the \emph{realizability problem} that is the moment problem on this space. Beside the standard inductive topology on the space of test functions $\C_c^\infty(\RR^d)$, we also represent this space as the uncountable intersection of weighted Sobolev spaces $H_k$ and we equip it with the associated strictly weaker projective topology. The corresponding space of generalized functions $\mathscr{D}'_{proj}(\RR^d)$ is strictly smaller than $\mathscr{D}'_{ind}(\RR^d)$, as it contains only generalized functions of finite order. The projective description is needed to apply the results presented in Subsection~\ref{Sec-MPNuclSpac} to the moment problem on $\mathcal{S}=\mathscr{D}'_{proj}(\RR^d)$. 

In Section~\ref{Sec-Main}, we formulate the main result of this paper, i.e.\! Theorem~\ref{MainThm1}. The only regularity assumption in the sense of Condition~\ref{condreg} is that the putative moments are generalized functions (in the projective topology, see Remark~\ref{RemD'Ind}). %Note that this requirement is equivalent to assuming that for each $n\in\NN$, the $n-$th moment function lies in the $n-$fold tensor product of the dual of one $H_k$, where the choice of the space may be different for each moment function. 
Furthermore, we assume a growth condition on the sequence of putative moment functions that expresses the conflicting nature of the Condition type~\ref{condquasi},~\ref{condreg} and~\ref{condgro} (see~Remark~\ref{RemarkE}). %~\ref{RemarkYuri}

In Subsection~\ref{Sec-MPNuclSpac}, we state the moment problem on the dual~$\Omega'$ of a general nuclear space $\Omega$ that is the projective limit of a family of separable Hilbert spaces and on subsets $\mathcal{S}$ of~$\Omega'$. We also recall the general result obtained by Y.~M.~Berezansky, Y.~G. Kondratiev and S.~N.~\v{S}ifrin for the moment problem on $\mathcal{S}=\Omega'$. We actually introduce their result under a slightly more general growth condition, which is given in Definition~\ref{DefSeq} (see~Remark~\ref{RemarkYuri}). This modification is essential to get the main result of this paper.
In Subsection~\ref{ProofMain}, we provide the detailed proof of Theorem~\ref{MainThm1}. Note that the theorem holds for the whole class of basic semi-algebraic sets of $\mathscr{D}'_{proj}(\RR^d)$, including the ones defined by an uncountable family of polynomials. To consider these kinds of sets, the inductive topology on $\C_c^\infty(\RR^d)$ plays an essential role, since $\mathcal{S}$ is closed, and so measurable, w.r.t.\!\! the strong topology on $\mathscr{D}'_{ind}(\RR^d)$ and the latter space is Radon (see Subsection~\ref{MeasD'}).

In Section \ref{Sec-Appl}, we use our main theorem to derive realizability results in more concrete cases. Fundamentally, given a specific desired support~$\mathcal{S}$, one has to find a representation of $\mathcal{S}$ as a basic semi-algebraic set of the space of generalized functions. Note that the result may depend on the chosen representation of $\mathcal{S}$. In Subsection~\ref{FinDimCase}, we describe how the new ideas employed in the proof of our main result allow us to extend the previous finite dimensional results to basic semi-algebraic sets of $\RR^d$ defined by an uncountable family of polynomials and to the most general bound of type \ref{condquasi}. In Subsection~\ref{1Subsec-Appl}, a more explicit description of the determinacy condition in terms of the scale of Sobolev spaces is introduced in the case when all moment functions are Radon measures. To avoid an extra unnecessary factorial factor in the determinacy bound obtained via Sobolev embedding (see Proposition~\ref{Prop1-1} and Remark \ref{RemarkNorm}), it is indispensable to use our more general definition of determining sequence which does not involve the norm of the moment functions as elements of the tensor product of the duals of the weighted Sobolev spaces. In Subsection~\ref{ExRadonMeas}, we investigate conditions under which such moment functions are realized by a random measure, that is by a finite measure concentrated on Radon measures. A spectral theoretical result of S.~N.~\v{S}ifrin \cite{S74} allows us also to weaken the determinacy condition. In Subsection~\ref{ExDensity} we show how to characterize, via moments, measures that are supported on the set of Radon measures with Radon-Nikodym density w.r.t.\! the Lebesgue measure fulfilling an \emph{a priori} $L^\infty$ bound. These examples also demonstrate that, in contrast to the finite dimensional case, a semi-algebraic set defined by uncountably many polynomials leads to very natural and treatable conditions on the moments in the infinite dimensional context. These positivity conditions can be seen as natural extensions of the classical conditions in the finite dimensional case, see Remarks~\ref{RemarkSDefCondMomMeas} and \ref{RemarkSDefCondMomMeas2}. In a forthcoming paper, we will treat further applications that require new additional ideas.

In Appendix~\ref{Sec-App1} and Appendix~\ref{Sec-App2}, we present some results from the theory of quasi-analyticity used in this paper and some considerations complementary to Subsection \ref{Sec-GenFunct}, respectively. Finally, in Appendix~\ref{Sec-App3} we give an explicit construction of a total subset of test functions fulfilling the requirement of the aforementioned determinacy condition. This construction allows us to obtain improved determinacy conditions in the particular cases considered in Section~\ref{Sec-Appl}.\\
We are convinced that the results contained in this paper are just the template for a multitude of forthcoming applications guided by their practical usefulness.

\section{Preliminaries}\label{Sec-Prel}
\subsection{The space of generalized functions}\label{Sec-GenFunct}\ \\
Let us first recall some standard general notations.

For $Y\subseteq \RR^d$ let us denote by $\mathcal{B}(Y)$ the Borel $\sigma$-algebra on $Y$, by $\C_c(Y)$ the space of all real-valued continuous functions on $\RR^d$ with compact support contained in $Y$ and by $\C_c^\infty(Y)$ its subspace of all infinitely differentiable functions. Moreover, $\C^+_c(Y)$ and $\C_c^{+,\infty}(Y)$ will denote the cones consisting of all non-negative functions in $\C_c(Y)$ and $\C_c^\infty(Y)$, respectively. Let $\NN_0$ be the set of all non-negative integers. For any
$\r=(r_1,\ldots,r_d)\in\RR^d$ and $\alpha=(\alpha_1,\ldots,\alpha_d)\in\NN_0^d$ one defines $\r^\alpha:=r_1^{\alpha_1}\cdots r_d^{\alpha_d}$. For any $\beta\in\NN_0^d$, the symbol $D^\beta$ denotes the partial derivative $\frac{\partial^{\left|\beta\right|}}{\partial r_1^{\beta_1}\cdots\partial r_d^{\beta_d}}$ where $\left|\beta\right|:=\sum_{i=1}^d\beta_i$. We will denote by $\Omega_\tau$ the space $\Omega$ endowed with a topology~$\tau$ and by $\Omega'_\tau$ its topological dual space. %The suffix will be dropped whenever the topology is clear from the context. 

In the following we introduce two different topologies on $\C_c^\infty(\RR^d)$, both making this space into a complete locally convex nuclear vector space.  

%\subsubsection{Topological structures on $\C_c^\infty(\RR^d)$}\label{Sec-GenFunct}\ \\
The classical topology considered on $\C_c^\infty(\RR^d)$ is the inductive topology $\tau_{ind}$, given by the standard construction of this space as the inductive limit of spaces of smooth functions with supports lying in an increasing sequence of compact subsets of $\RR^d$ (see~Definition~\ref{DefInd}). We denote by $\mathscr{D}_{ind}(\RR^d)$ the space $\C_c^\infty(\RR^d)$ equipped with~$\tau_{ind}$. On the other hand, the space $\C_c^\infty(\RR^d)$ can be also endowed with a projective topology~$\tau_{proj}$ in the following way (see~Definition~\ref{DProj} for an equivalent definition of~$\tau_{proj}$ and see \cite[Chapter I, Section 3.10]{B86} for more details). 
\begin{definition}\label{TeoBerez}\ \\
Let $I$ be the set of all $k=(k_1, k_2(\r))$ such that $k_1\in\NN_0$, $k_2\in\C^\infty(\RR^d)$ with $k_2(\r)\geq 1$ for all $\r\in\RR^d$.  For each $k=(k_1, k_2(\r))\in I$, consider the weighted Sobolev space $W_2^{k_1}(\RR^d, k_2(\r) d\r)$ defined as the completion of $\C_c^\infty(\RR^d)$ w.r.t.\! the following weighted norm
\begin{equation}\label{NormW2weighted}
\|\varphi\|_{W_2^{k_1}(\RR^d, k_2(\r) d\r)}:=\left(\sum_{|\beta|\leq k_1}\int_{\RR^d}\left |(D^{\beta}\varphi)(\r)\right |^2k_2(\r)d\r\right)^{\frac 12}.
\end{equation}
Then we define $$\mathscr{D}_{proj}(\RR^d):=\projlim\limits_{(k_1, k_2(\r))\in I}W_2^{k_1}(\RR^d, k_2(\r) d\r),
$$ and we denote by $\tau_{proj}$ the projective limit topology induced on $\C_c^\infty(\RR^d)$ by this construction.
\end{definition}

The previous definition of $\mathscr{D}_{proj}(\RR^d)$ is due to Y. M. Berezansky who also proved that such a projective limit is \emph{nuclear} (see \cite[Theorem~3.9, p.78]{B86} for the proof of this result). The latter property, as well as the construction of $\mathscr{D}_{proj}(\RR^d)$ as the projective limit of Hilbert spaces, is needed to apply the results of Subsection~\ref{Sec-MPNuclSpac} to the realizability problem on $\mathscr{D}'_{proj}(\RR^d)$.
 
Note that as sets, $\mathscr{D}_{ind}(\RR^d)$ and $\mathscr{D}_{proj}(\RR^d)$ coincide but the topologies $\tau_{ind}$ and $\tau_{proj}$ are not equivalent. In fact, it easily follows from the definitions of the two topologies that $\tau_{proj}\subset\tau_{ind}.$
Hence, we have that $\mathscr{D}'_{proj}(\RR^d)\subseteq\mathscr{D}'_{ind}(\RR^d)$ but this inclusion is actually strict. \\%For instance, let us consider in the case $d=1$ the functional $\eta(\varphi):=\sum_{n\in\NN_0} D^n\varphi(n)$ for any $\varphi\in\C_c^\infty(\RR)$. It is possible to prove that $\eta$ is in $\mathscr{D}'_{ind}(\RR)$ but not in $\mathscr{D}'_{proj}(\RR)$.\\

In what follows, $\mathscr{D}(\RR^d)$ and $\mathscr{D}'(\RR^d)$ are understood to be $\mathscr{D}_{proj}(\RR^d)$ and $\mathscr{D}'_{proj}(\RR^d)$, respectively. The suffix will be specified only whenever there might be ambiguity. We also denote by $\langle f, \eta\rangle$ the duality pairing between $\eta\in\mathscr{D}'(\RR^d)$ and $f\in\mathscr{D}(\RR^d)$ (see~\cite{B86, BeKo88} for more details). 

\subsection{Realizability problem on $\mathscr{D}'(\RR^d)$}\label{Subsec-RPonDPrime}\ \\
Let us introduce the main objects involved in the realizability problem on $\mathscr{D}'(\RR^d)$.

A generalized process on~$\mathscr{D}'(\RR^d)$ is a finite measure $\mu$ defined on the Borel $\sigma-$algebra on~$\mathscr{D}'(\RR^d)$. Moreover, we say that a generalized process $\mu$ is \emph{concentrated on} a measurable subset $\mathcal{S}\subseteq\mathscr{D}'(\RR^d)$ if $\mu\left(\mathscr{D}'(\RR^d)\setminus\mathcal{S}\right)=0$. 

\begin{definition}[Finite $n-$th local moment]\ \\
Given $n\in\NN$, a generalized process $\mu$ on $\mathscr{D}'(\RR^d)$ has \emph{finite $n-$th local moment} (or local moment of order $n$) if for every $f\in\C_c^\infty(\RR^d)$ we have
$$\int_{\mathscr{D}'(\RR^d)}|\langle f, \eta\rangle|^n \mu(d\eta)<\infty.$$
\end{definition}
%The latter condition is equivalent to the fact that 
%\begin{equation}\label{GenMomFunct}
% (f_1,\ldots, f_n)\mapsto\int_{\mathscr{D}'(\RR^d)} \langle f_1\otimes\cdots\otimes f_n, \eta^{\otimes n} \rangle \mu(d\eta).
% \end{equation}  
%is a well-defined multilinear functional on $\left(\C_c^\infty(\RR^d)\right)^{\times n}$. In fact, since $\mu$ has finite $n-$th local moment, for any $f_1,\ldots, f_n\in\C_c^\infty(\RR^d)$ we get
%$$
%\int_{\mathscr{D}'(\RR^d)} \langle f_1\otimes\cdots\otimes f_n, \eta^{\otimes n} \rangle \mu(d\eta)\leq\int_{\mathscr{D}'(\RR^d)} \prod_{i=1}^n |\langle f_i, \eta \rangle| \mu(d\eta)\leq \prod_{i=1}^n\left(\int_{\mathscr{D}'(\RR^d)}|\langle f_i, \eta \rangle|^n \mu(d\eta)\right)^{\frac 1n}\! \!\!<\!\infty.$$
%The functional in \eqref{GenMomFunct} is the \emph{$n-$th moment function} of $\mu$. In the following, we require slightly more regularity on the moment functions, but this assumption is easy to check in most of applications.

\begin{definition}[$n-$th generalized moment function]\ \\
Given $n\in\NN$, a generalized process $\mu$ on $\mathscr{D}'(\RR^d)$ has $n-$th generalized moment function in the sense of $\mathscr{D}'(\RR^d)$ if $\mu$ has finite $n-$th local moment and if the functional $f\mapsto\int_{\mathscr{D}'(\RR^d)}|\langle f, \eta\rangle|^n \mu(d\eta)$ is continuous on $\mathscr{D}(\RR^d)$.\\ In fact, by the Kernel Theorem, for such a generalized process $\mu$ there exists a symmetric functional $m^{(n)}_{\mu}\in\mathscr{D}'(\RR^{dn})$, which will be called the \emph{$n-$th generalized moment function in the sense of $\mathscr{D}'(\RR^d)$}, such that
%for any $f_1,\ldots,f_n\in\C_c^\infty(\RR^d)$ the following holds 
%\begin{equation*}\label{nMomTensor}
% \langle f_1\otimes\cdots\otimes f_n, m_{\mu}^{(n)} \rangle =\int_{\mathscr{D}'(\RR^d)} \langle f_1\otimes\cdots\otimes f_n, \eta^{\otimes n} \rangle \mu(d\eta).
% \end{equation*} 
for any $f^{(n)}\in\C_c^\infty(\RR^{dn})$ we have
\begin{equation}\label{nMomTensor}
\langle f^{(n)}, m_{\mu}^{(n)} \rangle =\int_{\mathscr{D}'(\RR^d)} \langle f^{(n)}, \eta^{\otimes n} \rangle \mu(d\eta).
 \end{equation} 
By convention, $m_{\mu}^{(0)}:=\mu(\mathscr{D}'(\RR^d))$.
\end{definition}
 
%\begin{prop}\ \\
%If $\mu$ is a generalized process on $\mathscr{D}'(\RR^d)$ with generalized moment functions (in the sense of $\mathscr{D}'(\RR^d)$) of any order, then for any $n\in\NN$ and for any $f^{(n)}\in(\left(\C_c^\infty(\RR^d)\right)^{\otimes n}$ we have
%\begin{equation*}
%\int_{\mathscr{D}'(\RR^d)}| \langle f^{(n)}, \eta^{\otimes n} \rangle| \mu(d\eta)<\infty\qquad\text{and}\qquad\langle f^{(n)}, m_{\mu}^{(n)} \rangle =\int_{\mathscr{D}'(\RR^d)} \langle f^{(n)}, \eta^{\otimes n} \rangle \mu(d\eta).
%\end{equation*}
%\end{prop}

For a generalized processes $\mu$ the generalized moment functions~$m_{\mu}^{(n)}$ are given by \eqref{nMomTensor}. The moment problem, which in an infinite dimensional context is often called the \emph{realizability problem}, addresses exactly the inverse question.
\begin{probl}[Realizability problem on $\mathcal{S}\subseteq\mathscr{D}'(\RR^d)$]\label{RealProb}\  \\
Let $\mathcal{S}$ be a measurable subset of $\mathscr{D}'(\RR^d)$, $N\in\mathbb{N}_0\cup\{+\infty\}$ and $m=(m^{(n)})_{n=0}^N$ such that each $m^{(n)}\in\mathscr{D}'(\RR^{dn})$ is a symmetric functional. Find a generalized process $\mu$ with generalized moments (in the sense of $\mathscr{D}'(\RR^d)$) of any order and concentrated on~$\mathcal{S}$ such that $$m^{(n)}=m^{(n)}_\mu\quad\text{for}\,\,\, n=0,\ldots, N,$$
i.e. $m^{(n)}$ is the $n-$th generalized moment function of $\mu$ for $n=0,\ldots, N$.
\end{probl}
If such a measure $\mu$ does exist we say that $(m^{(n)})_{n=0}^N$ is \emph{realized} by $\mu$ on~$\mathcal{S}$. Note that the definition requires that one finds a measure concentrated on~$\mathcal{S}$ and not only on $\mathscr{D}'(\RR^d)$. In the case $N=\infty$ one speaks of the ``full realizability problem'', otherwise of the ``truncated realizability problem''.

\section{Realizability problem on basic semi-algebraic subsets of $\mathscr{D}'(\RR^d)$}\label{Sec-Main}
To simplify the notation in the following we denote by
$\mathcal{M}^*(\mathcal{S})$ the collection of all generalized processes concentrated on a measurable subset $\mathcal{S}$~of~$\mathscr{D}'(\RR^d)$ with generalized moment functions (in the sense of $\mathscr{D}'(\RR^d)$) of any order and by $\mathcal{F}\left(\mathscr{D}'(\RR^d)\right)$ the collection of all infinite sequences $(m^{(n)})_{n\in\NN_0}$ such that each $m^{(n)}\in\mathscr{D}'(\RR^{dn})$ is a symmetric functional of its $n$ variables.%, namely the tensor product $\left(\Omega'\right)^{\otimes n}$ is considered to be symmetric.

Let $\mathscr{P}_{\C^\infty_c}\left(\mathscr{D}'(\RR^d)\right)$ be the set of all polynomials on $\mathscr{D}'(\RR^d)$ of the form 
\begin{equation}\label{polyDprime}
P(\eta) := \sum_{j=0}^N\langle p^{(j)},\eta^{\otimes j}\rangle,
\end{equation}
where $p^{(0)}\in\RR$ and $p^{(j)}\in\C^\infty_c(\RR^{dj})$, $j=1,\ldots,N$ with $N\in\NN$. 

A subset $\mathcal{S}$ of $\mathscr{D}'(\RR^d)$ is said to be \emph{basic semi-algebraic} if it can be written as
\begin{equation}\label{S-Semialg}
\mathcal{S}=\bigcap_{i\in Y}\left\{\eta\in\mathscr{D}'(\RR^d)| \ P_i(\eta)\geq 0\right\},
\end{equation} 
where $Y$ is an index set and $P_i\in\mathscr{P}_{\C^\infty_c}\left(\mathscr{D}'(\RR^d)\right)$.  Note that the index set $Y$ is not necessarily countable. Moreover, let $\mathscr{P}_{\mathcal{S}}$ be the set of all the polynomials $P_i$'s defining $\mathcal{S}$. W.l.o.g. we assume that $P_0$ is the constant polynomial $P_0(\eta)=1$ for all $\eta\in\mathscr{D}'(\RR^d)$ and that~$0 \in Y$. \\

In the following, we are going to investigate the full realizability problem (see~Problem~\ref{RealProb}) on $\mathcal{S}$ of the form~\eqref{S-Semialg}.\\

First let us introduce the concept of \emph{determining} sequence, which essentially is a growth condition on the sequence of the $m^{(n)}$'s.  We will see that this property gives the uniqueness of the realizing measure.

\begin{definition}[Determining sequence]\label{DefSeq}\ \\
Let $m\in\mathcal{F}\left(\mathscr{D}'(\RR^{d})\right)$ and $E$ be a total subset of $\mathscr{D}(\RR^{d})$, i.e.\! the linear span of $E$ is dense in~$\mathscr{D}(\RR^{d})$. Let us define the sequence $(m_n)_{n\in\NN_0}$ as follows
\begin{equation}\label{defCond}
m_0:=\sqrt{|m^{(0)}|}\,\text{ and }\,m_n:= \sqrt{\sup_{f_1,\ldots,f_{2n}\in E}|\langle f_1\otimes\cdots\otimes f_{2n},m^{(2n)}\rangle|},\, \forall\,n\geq 1.
\end{equation}
The sequence $m$ is said to be \emph{determining} if and only if there exists a total subset $E$ of $\C_c^\infty(\RR^d)$ such that for any $n\in\NN_0$, $m_n<\infty$ and the class $C\{m_n\}$ is quasi-analytic (see\!~Definition~\ref{QuasianalityCLA} and Theorem~\ref{DJ-C}).
\end{definition}

The version of the Riesz linear functional for the moment problem on $\mathscr{D}_{proj}'(\RR^d)$ is given by the following. 
\begin{definition}\label{DefFunct}\ \\
Given $m\in\mathcal{F}\left(\mathscr{D}'(\RR^{d})\right)$, we define its associated Riesz functional $L_m$ as 
\begin{eqnarray*}
L_m: &\mathscr{P}_{\C_c^\infty}\left(\mathscr{D}'(\RR^d)\right)&\to\RR \\
 \  &P(\eta)=\sum\limits_{n=0}^N\langle p^{(n)} , \eta^{\otimes n}\rangle & \mapsto L_m(P):=\sum_{n=0}^N\langle p^{(n)} , m^{(n)}\rangle .
\end{eqnarray*}
\end{definition}
Note that in the case when the sequence $m$ is realized by a non-negative measure $\mu\in\mathcal{M}^*(\mathcal{S})$ on a subset $\mathcal{S}\subseteq\mathscr{D}'(\RR^d)$, a direct calculation shows that for any polynomial $P\in\mathscr{P}_{\C_c^\infty}(\mathscr{D}'(\RR^d))$ 
\begin{equation}
L_m(P)=\int_{\mathcal{S}}{P(\eta)\,\mu(d\eta)}.\label{intFormFunct}
\end{equation}

The Riesz functional allows us to state our main result in a concise form.
\begin{theorem}\label{MainThm1}\ \\
Let $m\in\mathcal{F}\left(\mathscr{D}'(\RR^{d})\right)$ be determining and $\mathcal{S}$ be a basic semi-algebraic set of the form~\eqref{S-Semialg}.
Then $m$ is realized by a unique non-negative measure $\mu\in\mathcal{M}^*(\mathcal{S})$ if and only if the following inequalities hold
\begin{equation}\label{MainThmCond1-1}
L_m(h^2)\geq 0, \,\,L_m(P_i h^2)\geq 0\,\,,\,\, \forall h\in\mathscr{P}_{\C_c^\infty}\left(\mathscr{D}'(\RR^d)\right),\,\forall i\in Y.
\end{equation}
\end{theorem}

In other words one can see the solution to the realizability problem as a way to read off from the moment functions support properties for any realizing measure.

\begin{remark}\ \\
Condition \eqref{MainThmCond1-1} is equivalent to require that the functional $L_m$ is non-negative on the quadratic module $\mathcal{Q}(\mathscr{P}_{\mathcal{S}})$.
We define the quadratic module $\mathcal{Q}(\mathscr{P}_{\mathcal{S}})$ associated to the representation \eqref{S-Semialg} of $\mathcal{S}$ as the convex cone in $\mathscr{P}_{\C_c^\infty}(\mathscr{D}'(\RR^d))$ given by
$$\mathcal{Q}(\mathscr{P}_{\mathcal{S}}):=\bigcup_{\stackrel{Y_0\subset Y}{|Y_0|<\infty}}\left\{\sum_{i\in Y_0} Q_i P_i\,:\, Q_i\in\Sigma_{\C_c^\infty}(\mathscr{D}'(\RR^d))\right\},$$
where $\Sigma_{\C^\infty_c}(\mathscr{D}'(\RR^d))$ denotes the subset of all polynomials in $\mathscr{P}_{\C^\infty_c}\left(\mathscr{D}'(\RR^d)\right)$ which can be written as sum of squares of polynomials.
\end{remark}

The determinacy condition given in Definition \ref{DefSeq} seems to be abstract, but it becomes actually very concrete whenever one can explicitly construct the set $E$. In fact, for any $n\in\NN_0$ we have that
\begin{equation}\label{tensordualD}
\mathscr{D}'(\RR^{dn})=\bigcup_{k\in I}(H'_{k})^{\otimes n}=\bigcup_{k\in I}H_{-k}^{\otimes n},
\end{equation}
where $k=(k_1,k_2(\r))\in I$, $H_{k}:=W_2^{k_1}(\RR^d, k_2(\r) d\r)$ and $I$ is as in Definition~\ref{TeoBerez} (for $n=0$, $H_{-k}^{\otimes n}=\RR$). From \eqref{tensordualD} follows that for any sequence $m\in\mathcal{F}\left(\mathscr{D}'(\RR^{d})\right)$ there exists a sequence $(k^{(n)})_{n\in\NN_0}\subset I$ such that for any $n\in\NN_0$ we get $m^{(n)}\in H_{-k^{(n)}}^{\otimes n}$. If we denote by $d(k^{(n)}, E):=\sup\limits_{f\in E}\|f\|_{H_{k^{(n)}}}$, then for the $m_n$'s defined in \eqref{defCond} we have
\begin{equation}\label{boundDeterm}
m_n\leq (d(k^{(2n)}, E))^{n}\|m^{(2n)}\|_{H_{-k^{(2n)}}^{\otimes{2n}}}^{\frac 12}.
\end{equation}
Hence, we can see that a preferable choice for $E$ is the one for which $\left(d(k^{(2n)},E)\right)_{n\in\NN}$ grows as little as possible. Such an $E$ can be obtained by using the following result, which we proved with a technique similar to the one of \cite[Chapter 4, Section 9]{Gel-Shi68} (see Appendix~\ref{Sec-App3} for the proof of Lemma~\ref{LemmaE}). 

\begin{lemma}\label{LemmaE} \ \\
Let $(c_n)_{n\in\NN_0}$ be an increasing sequence of positive numbers which is not quasi-analytic and let $m\in\mathcal{F}( \mathscr{D}'(\mathbb{R}^d) )$. For any $n\in\NN_0$, let $k^{(n)}:=(k_1^{(n)}, k_2^{(n)})\in I$ be such that $m^{(n)}\in H_{-k^{(n)}}^{\otimes n}$ where $H_{k^{(n)}}:=W_2^{k_1^{(n)}}\!\!\!(\RR^d, k_2^{(n)}(\r) d\r)$ and $I$ is as in Definition~\ref{TeoBerez}. 
Then the set
\[
E := \left\{f \in \mathscr{D}(\mathbb{R}^d) \left| \forall\ n\in\NN_0,\  \| f\|_{H_{k^{(n)}}} \leq c^d_{k_1^{(n)}}  \sup_{\stackrel{\z\in\RR^d}{\|\z\|\leq n}}\sup_{\x \in [-1,1]^d}\sqrt{ k_2^{(n)}(\z+\x)} \right.\right\}
\]
is total in $\mathscr{D}(\mathbb{R}^d)$.
\end{lemma}
For such a set $E$, using \eqref{boundDeterm}, we get that $$m_n\leq c^{dn}_{k_1^{(n)}}  \left(\sup_{\stackrel{\z\in\RR^d}{\|\z\|\leq n}}\sup_{\x \in [-1,1]^d}\sqrt{ k_2^{(n)}(\z+\x)}\right)^n\|m^{(2n)}\|_{H_{-k^{(2n)}}^{\otimes{2n}}}^{\frac 12}.$$
Note that concrete examples of increasing sequences of positive numbers which are not quasi-analytic are provided in Appendix~\ref{Sec-App3}. 

\begin{remark}\label{RemarkE}\ \\
The more regularity is known on the sequence $m$ the weaker is the restriction on the growth of the $m^{(2n)}$ required in Theorem \ref{MainThm1}. Let us discuss two extremal cases.
\begin{itemize}
\item If each $m^{(n)}$ is in $H_{-k}^{\otimes n}$ where $k=(k_1, k_2(\r))\in I$ with both $k_1$ and $k_2$ independent of $n$, then both $c_{k_1^{(n)}}$ and $\sup\limits_{\stackrel{\z\in\RR^d}{\|\z\|\leq n}}\sup\limits_{\x \in [-1,1]^d}\sqrt{ k_2^{(n)}(\z+\x)}$ in Lemma~\ref{LemmaE} are constant w.r.t.\! $n$ and so a sufficient condition for the determinacy of $m$ is the quasi-analyticity of the class $C\{\left\|m^{(2n)}\right\|_{H_{-k}^{\otimes 2n}}^{1/2}\}$.
\item If each $m^{(n)}$ is in $H_{-k^{(n)}}^{\otimes n}$ where $k^{(n)}=(k_1, k_2^{(n)}(\r))\in I$ with $k_1$ independent of $n$, then $c_{k_1^{(n)}}$ in Lemma~\ref{LemmaE} is constant w.r.t.\! $n$ and so a sufficient condition for the determinacy of $m$ is the quasi-analyticity of the class $$C\left\{\left(\sup_{\stackrel{\z\in\RR^d}{\|\z\|\leq n}}\sup_{\x \in [-1,1]^d}\sqrt{ k_2^{(n)}(\z+\x)}\right)^n\left\|m^{(2n)}\right\|_{H_{-k^{(n)}}^{\otimes 2n}}^{1/2}\right\}.$$
Hence, the condition on $m$ of being determining also contains the growth of the sequence of functions $(k_2^{(n)})_{n\in\NN}$. For a concrete application of this to moment functions which are themselves Radon measures see Subsection~\ref{1Subsec-Appl}.
\end{itemize}
\end{remark}

\section{Proof of the main result}\label{ProofMainResult}\ 
The proof of the main result of this paper, Theorem~\ref{MainThm1}, is based on the application of a general theorem about the realizability problem on nuclear spaces to $\mathscr{D}'_{proj}(\RR^d)$. Such a result guarantees, under our assumptions on the starting sequence $m$, the existence and the uniqueness of a realizing measure on $\mathscr{D}'_{proj}(\RR^d)$ (see Subsection~\ref{Sec-MPNuclSpac}). The main point of Theorem~\ref{MainThm1} is to show that the realizing measure is actually supported on $\mathcal{S}$. In the case when the pre-given semi-algebraic set $\mathcal{S}$ is defined by an uncountable family of polynomials, we need to consider the inductive topology on $\C_c^\infty(\RR^d)$ in order to prove the support properties. The inductive topology plays an essential role in this case, since $\mathcal{S}$ is closed w.r.t.\! the strong topology on $\mathscr{D}'_{ind}(\RR^d)$ and the latter space is Radon (see Subsection~\ref{MeasD'}).

Before giving the proof of the main theorem, we need to describe the general framework in more details and give some preliminary results.

\subsection{Realizability problem on nuclear spaces}\label{Sec-MPNuclSpac}\ \\
In the following we will consider all the spaces as being separable and real.

Let us consider a family $(H_k)_{k\in K}$ of Hilbert spaces ($K$ is an index set containing~$0$) which is directed by topological embedding, i.e. 
$$\forall\ k_1,k_2\in K\,\,\exists\, k_3\,:\, H_{k_3}\subseteq H_{k_1}\, ,\, H_{k_3}\subseteq H_{k_2}.$$
We assume that each $H_k$ is embedded topologically into $H_0$. Note that the $H_k$'s are not necessarily Sobolev spaces.

Let $\Omega$ be the projective limit of the family $(H_k)_{k\in K}$ endowed with the associated projective limit topology and let us assume that $\Omega$ is nuclear, i.e.\! for each $k_1\in K$ there exists $k_2\in K$ such that the embedding $H_{k_2}\subseteq H_{k_1}$ is quasi-nuclear.

Let us denote by $\Omega'$ the topological dual space of $\Omega$. We control the classical rigging by identifying $H_0$ and its dual $H'_0$. With this identification one can define the duality pairing between elements in $H_k$ and in its dual $H'_k=H_{-k}$ using the inner product in $H_0$. For this reason, in the following we will denote by $\langle f, \eta\rangle$ the duality pairing between $\eta\in\Omega'$ and $f\in\Omega$ (see~\cite{B86, BeKo88} for more details).

Consider the $n-$th ($n\in\NN_0$) tensor power $\Omega^{\otimes n}$ of the space $\Omega$ which is defined as the projective limit of $H_k^{\otimes n}$; for $n=0$, $H_k^{\otimes n}=\RR$. Then its dual space is
\begin{equation}\label{tensorDual}
\left(\Omega^{\otimes n}\right)'=\bigcup_{k\in K}\left(H_k^{\otimes n}\right)'=\bigcup_{k\in K}(H'_k)^{\otimes n}=\bigcup_{k\in K}H_{-k}^{\otimes n},
\end{equation}
which we can equip with the weak topology. 
%Assume that an involution $\Omega\ni f\mapsto \overline{f}\ni\Omega$ is defined in $\Omega$ and that it can be extended by continuity to an involution in each $H_k$ for all $k\in K$ and in a natural way to tensor powers and dual spaces.\\
%In the following, we will consider all our spaces real, namely all the elements are invariant w.r.t. the involution, i.e. $f=\overline{f}$.

All the definitions about the realizability problem on $\mathscr{D}'(\RR^d)$ given in Subsection~\ref{Subsec-RPonDPrime} can be straightforwardly generalized by replacing $\mathscr{D}(\RR^d)$ by $\Omega$ and $\mathscr{D}(\RR^{dn})$ by $\Omega^{\otimes n}$. In this more general context, $\mathcal{F}(\Omega')$ denotes the collection of all infinite sequences $(m^{(n)})_{n\in\NN_0}$ such that each $m^{(n)}\in\left(\Omega^{\otimes n}\right)'$ is a symmetric functional, namely an element of the symmetric $n-$fold tensor product of $\Omega'$.

An obvious positivity property which is necessary for an element in $\mathcal{F}(\Omega')$ to be the moment sequence of some measure on $\Omega'$ is the following.
\begin{definition}[Positive semidefinite sequence]\label{PosSemiDef}\ \\
A sequence $m\in\mathcal{F}(\Omega')$ is said to be \emph{positive semidefinite} if for any $f^{(j)}\in\Omega^{\otimes j}$ we have
$$\sum_{j,l=0}^\infty \langle  f^{(j)} \otimes f^{(l)}, m^{(j+l)}\rangle\geq 0.$$
\end{definition}
This is a straightforward generalization of the classical notion of positive semidefiniteness of the Hankel matrices considered in the finite dimensional moment problem, that is equivalent to require that the associated Riesz functional is non-negative on squares of polynomials. Note that, as we work with real spaces, we choose the involution on $\Omega$ considered in~\cite{BeKo88} to be the identity.

The definition of determining sequence is the obvious analogous of Definition~\ref{DefSeq} for a sequence $m\in\mathcal{F}(\Omega')$ and so, using \eqref{tensorDual}, we get \eqref{boundDeterm} in this general case. However, the explicit construction of a subset $E$ of $\Omega$ for which $\left(d(k^{(2n)},E)\right)_{n\in\NN}$ grows as little as possible must depend on the structure of the Hilbert spaces $H_k$. Hence, an analogous construction to the one in Lemma~\ref{LemmaE} cannot be given in abstract but it will always depend on the concrete structure of the particular $H_k$'s.\\

Let us state now the fundamental result for the full realizability problem in the case $\mathcal{S}=\Omega'$ and $\Omega'$ is a Suslin space (see\!~\cite[Vol.~II, Theorem~2.1, p.54]{BeKo88} and \cite{BS71}).
\begin{theorem}\label{KondrThm}\ \\
If $m\in\mathcal{F}(\Omega')$ is a positive semidefinite sequence which is also determining, then there exists a unique non-negative generalized process $\mu\in\mathcal{M}^*(\Omega')$ such that for any $f^{(n)}\in\Omega^{\otimes n}$
 \begin{equation*}\label{GenRealiz}
 \left\langle f^{(n)}, m^{(n)} \right\rangle=\int_{\Omega'}\left\langle f^{(n)},\eta^{\otimes n}\right\rangle\mu( d\eta).
 \end{equation*}
 \end{theorem}

\begin{remark}\label{RemarkYuri}\ \\
The original proof of Theorem~\ref{KondrThm} in \cite{BeKo88} uses a slightly less general definition of determining sequence. Indeed, the authors require that the class $$C\left\{d(k^{(2n)},E)^n\left\|m^{(2n)}\right\|_{H_{-k^{(2n)}}^{\otimes 2n}}^{1/2}\right\}$$ is quasi-analytic, which in turn implies that the class $C\{m_n\}$ is also quasi-analytic. Nevertheless, their proof also applies just using the bound given by Definition~\ref{DefSeq} for $m\in\mathcal{F}(\Omega')$. The latter has actually the advantage to guarantee, whenever $m$ is realizable on $\Omega$, the log-convexity of the sequence $(m_n)_{n\in\NN_0}$. This property is essential in the proof of the main result of this paper.\\
Let us also note that the proof of Theorem~\ref{KondrThm} actually shows that the measure $\mu$ is concentrated on one of the Hilbert spaces $H_{-k'}$ for some index $k'\in K$ depending on the sequence~$m$. Indeed, the index $k'$ is the one such that the embedding of $H_{k'}$ into $H_{k^{(2)}}$ is quasi-nuclear (see \cite[Remark~1,~p.~72]{BeKo88}). However, note that the assumptions of Theorem~\ref{KondrThm} do not require that all $m^{(n)}\in H_{-k'}^{\otimes n}$.
\end{remark}

In the following we are going to apply Theorem~\ref{KondrThm} for $\Omega=\mathscr{D}(\RR^d)$ constructed as the projective limit of a family of weighted Sobolev spaces $H_k:=W_2^{k_1}(\RR^d, k_2(\r) d\r)$, which is nuclear (see Subsection~\ref{Sec-GenFunct}). %for any $k=(k_1, k_2(\r))\in I$
Since $\Omega^{\otimes n}=\mathscr{D}(\RR^{dn})$, in this case the sequence $m$ consists of symmetric generalized functions, i.e.\!~$m^{(n)}\in\mathscr{D}'(\RR^{dn})$. Theorem~\ref{KondrThm} gives a solution for the full realizability problem on $\mathcal{S}=\mathscr{D}'(\RR^d)$ whenever the sequence $m$ is positive semidefinite and determining.

\subsection{Measurability of $\mathscr{D}'_{proj}(\RR^d)$ in $\mathscr{D}'_{ind}(\RR^d)$}\label{MeasD'}\ \\
The weak topology $\tau^{proj}_w$ [$\tau^{ind}_w$, resp.] on $\mathscr{D}_{proj}'(\RR^d)$ [$\mathscr{D}'_{ind}(\RR^d)$, resp.] is the smallest topology such that the mappings
$\eta \mapsto \langle f, \eta \rangle$
are continuous for all $f \in \C^\infty_c(\RR^d)$. It is easy to see that ${\tau}_w^{proj}$ coincides with the relative topology %~$\breve{\tau}_w^{ind}$ 
given by $\tau_w^{ind}$ on~$\mathscr{D}'_{proj}(\RR^d)\subset\mathscr{D}'_{ind}(\RR^d)$. As a consequence, the Borel $\sigma-$algebras generated by these two topologies also coincide and we can easily conclude that
\begin{equation}\label{CorSigmaAlgebra}
\sigma({\tau}_w^{proj})=\sigma({\tau}_w^{ind})\cap\mathscr{D}'_{proj}(\RR^d).
\end{equation}

Let us recall some properties of $\mathscr{D}'_{ind}(\RR^d)$.\\
Consider the strong topology $\tau_s^{ind}$ on $\mathscr{D}'_{ind}(\RR^d)$. It is well known that $\tau_s^{ind}$ coincides with the topology of compact convergence $\tau_c^{ind}$ and so, by Corollary~1 in \cite[Chapter II,~p.115]{Sch73}, $\left(\mathscr{D}'_{ind}(\RR^d),\tau_c^{ind}\right)$ is Lusin. Moreover, since $\tau_w^{ind}\subset\tau_s^{ind}$, the space $\left(\mathscr{D}'_{ind}(\RR^d),\tau_w^{ind}\right)$ is also Lusin. Hence, by Theorem~9 in \cite[Chapter~II,~p.122]{Sch73}, the following proposition holds.
\begin{prop}\label{D'Radon}\ \\
$\left(\mathscr{D}'_{ind}(\RR^d),\tau_w^{ind}\right)$ is a Radon space, i.e.\! every finite Borel measure on $\mathscr{D}'_{ind}(\RR^d)$ is inner regular.
\end{prop}

We were unable to find in the literature an analogous result establishing whether $\left(\mathscr{D}'_{proj}(\RR^d),\tau_w^{proj}\right)$ is a Radon space or not. In fact, the techniques used in \cite{Sch73} do not apply to $\mathscr{D}'_{proj}(\RR^d)$. 

On the level of Borel $\sigma-$algebras on $\mathscr{D}'_{ind}(\RR^d)$, we have that any Borel $\sigma-$algebra generated by a topology weaker than $\tau_s^{ind}$ coincides with the one generated by $\tau_s^{ind}$, since $\left(\mathscr{D}'_{ind}(\RR^d), \tau_s^{ind}\right)$ is a Lusin space and so Suslin (see \cite[Corollary~2, p.101]{Sch73}).
 
\subsection{Proof of the main result}\label{ProofMain}\ \\
Let us first note that the definition of $\mathscr{P}_{\C^\infty_c}\left(\mathscr{D}'_{proj}(\RR^d)\right)$ given in \eqref{polyDprime} can be extended to $\mathscr{P}_{\C^\infty_c}\left(\mathscr{D}'_{ind}(\RR^d)\right)$ by taking $\eta\in\mathscr{D}'_{ind}(\RR^d)$.
\begin{prop}\ \\
Every polynomial in $\mathscr{P}_{\C^\infty_c}\left(\mathscr{D}'_{ind}(\RR^d)\right)$ is continuous w.r.t.\! $\tau_s^{ind}$.
Hence, the basic semi-algebraic set $\mathcal{S}$ defined in \eqref{S-Semialg} is closed in $(\mathscr{D}'_{ind}(\RR^d), \tau_s^{ind})$.
\end{prop}
\proof\ \\
To show the continuity of a generic polynomial of the form \eqref{polyDprime}, it suffices to prove that for all $j\in\NN$ the functions 
\begin{eqnarray*}
\mathscr{D}'_{ind}(\RR^d)&\to&\RR\nonumber \\
 \eta & \mapsto& \langle p^{(j)} , \eta^{\otimes j}\rangle\label{mappa}
\end{eqnarray*}
are continuous w.r.t. $\tau_s^{ind}$.

For any fixed $j\in\NN$, we first consider the mapping $\eta\mapsto\eta^{\otimes j}$ which is continuous as a function from the space $(\mathscr{D}'_{ind}(\RR^d),\tau_s^{ind})$ to the algebraic tensor product $\left(\mathscr{D}'_{ind}(\RR^d)\right)^{\otimes j}$ endowed with the $\pi-$topology (see \cite[Definition~43.2]{Tre67}). 
Moreover, the closure of the latter space is isomorphic to $(\mathscr{D}'_{ind}(\RR^{jd}), \tau_s^{ind})$ (see~\cite[Theorem~51.7]{Tre67}). Finally, the function 
$ \zeta \mapsto \langle p^{(j)} , \zeta\rangle$ on $\mathscr{D}_{ind}'(\RR^{jd})$ is continuous w.r.t.\! the weak topology on this space and hence, it is also continuous w.r.t.\! the strong one.\\
\endproof

\begin{corollary}\label{S-Meas}\ \\
The semi-algebraic set $\mathcal{S}$ defined as in \eqref{S-Semialg} is measurable w.r.t. the Borel $\sigma-$algebra  $\sigma(\tau_w^{ind})$ generated by the weak topology on $\mathscr{D}'_{ind}(\RR^d)$.
\end{corollary}
\proof\ \\
The previous proposition implies that $\mathcal{S}\in\sigma(\tau_s^{ind})$.  As $(\mathscr{D}'_{ind}(\RR^d), \tau_s^{ind})$ is a Lusin space and so Suslin, $\sigma(\tau_w^{ind})$ and $\sigma(\tau_s^{ind})$ coincide (see \cite[Corollary~2,~p.101]{Sch73}). Hence, $\mathcal{S}\in\sigma(\tau_w^{ind})$.\\
\endproof

Before proving Theorem \ref{MainThm1} we need to show some preliminary results. Remind that throughout the whole section we consider a sequence $m\in\mathcal{F}\left(\mathscr{D}'_{proj}(\RR^{d})\right)$.
\begin{definition}\label{SHIFT}\ \\
Given a polynomial $P\in\mathscr{P}_{\C_c^\infty}(\mathscr{D}'_{proj}(\RR^d))$ of the form $P(\eta) := \sum_{j=0}^{N}\langle p^{(j)},\eta^{\otimes j}\rangle$, we define the sequence $_Pm=\left((_Pm)^{(n)}\right)_{n\in\NN_0}$ in $\mathcal{F}\left(\mathscr{D}'_{proj}(\RR^{d})\right)$ as follows
$$
\forall f^{(n)}\in\C_c^\infty(\RR^{nd}),\quad\langle f^{(n)},(_P m)^{(n)}\rangle:=\sum_{j=0}^N\langle p^{(j)}\otimes f^{(n)}, m^{(n+j)}\rangle.
$$
\end{definition}
In terms of the Riesz functional introduced in Definition~\ref{DefFunct}, the previous definition takes the following form
\begin{equation}\label{shiftProperty}
\forall P,Q\in\mathscr{P}_{\C_c^\infty}(\mathscr{D}'_{proj}(\RR^d)),\quad L_{_{P}m}(Q):=L_m(PQ).
\end{equation}
\begin{remark}\label{semidefFunct}\ \\
The conditions \eqref{MainThmCond1-1} can be interpreted as that the sequence $(m^{(n)})_{n\in\NN_0}$ and all its shifted versions $((_{P_i}m)^{(n)})_{n\in\NN_0}$ are positive semidefinite in the sense of Definition~\ref{PosSemiDef}.
\end{remark}

\begin{lemma}\label{Prop2-1}\ \\
Let $P\in\mathscr{P}_{\C_c^\infty}(\mathscr{D}'_{proj}(\RR^d))$. If $m$ is realized on $\mathscr{D}_{proj}'(\RR^d)$ by a non-negative measure $\mu\in\mathcal{M}^*(\mathscr{D}_{proj}'(\RR^d))$, then the sequence $_P m$ is realized by the signed measure $P\mu$ on $\mathscr{D}_{proj}'(\RR^d)$.
\end{lemma}
\proof\ \\
Let $n\in\NN$ and $Q(\eta) :=\langle f^{(n)}, \eta^{\otimes n} \rangle$ with $f^{(n)}\in\C_c^\infty(\RR^{nd})$. Then, using~\eqref{intFormFunct} and~\eqref{shiftProperty}, one gets that 
$$
\int_{\mathscr{D}_{proj}'(\RR^d)} \langle f^{(n)}, \eta^{\otimes n} \rangle P(\eta)\mu(d\eta)= L_m(QP)=L_{_{P}m}(Q)=\langle f^{(n)}, (_P m)^{(n)} \rangle.
$$
%\begin{eqnarray*}
%\langle f^{(n)}, (_P m)^{(n)} \rangle&=&\sum_{j=0}^N \langle p^{(j)}\otimes f^{(n)}, m^{(n+j)}\rangle\\
%&=&\sum_{j=0}^N\int_{\mathscr{D}_{proj}'(\RR^d)} \langle p^{(j)} \otimes f^{(n)}, \eta^{\otimes (n+j)} \rangle \mu(d\eta)\\
%&=&\int_{\mathscr{D}_{proj}'(\RR^d)} \langle  f^{(n)}, \eta^{\otimes n} \rangle\sum_{j=0}^N \langle p^{(j)}, \eta^{\otimes j} \rangle \mu(d\eta)\\
%&=&\int_{\mathscr{D}_{proj}'(\RR^d)} \langle f^{(n)}, \eta^{\otimes n} \rangle P(\eta)\mu(d\eta).
%\end{eqnarray*}
\endproof

\begin{prop}\label{Prop3-1}\ \\
If $m$ is realized by a measure $\mu\in\mathcal{M}^*(\mathscr{D}_{proj}'(\RR^{d}))$ and $m$ is determining, then the sequence $_Pm$ is also determining. \end{prop}

\proof\ \\
Let us first recall that
$\mathscr{D}_{proj}(\RR^d)=\projlim\limits_{k\in I}H_k,$
where $I$ is as in Definition~\ref{TeoBerez} and $H_k:=W_2^{k_1}(\RR^d, k_2(\r) d\r)$ for any $k=(k_1, k_2(\r))\in I$ (see\!~Subsection~\ref{Sec-GenFunct}).\\

Since $m$ is determining in the sense of Definition~\ref{DefSeq}, there exists a subset $E$ total in $\mathscr{D}_{proj}(\RR^d)$ such that for any $n\in\NN_0$, $m_n<\infty$ and the class $C\{{m}_n\}$ is quasi-analytic, where
$$
{m}_n:=\sqrt{\sup_{f_1,\ldots,f_{2n}\in E}\left|\langle f_1\otimes\cdots\otimes f_{2n}, m^{(2n)}\rangle\right|}.
$$
It is easy to see that, since $m$ is realized by a measure $\mu\in\mathcal{M}^*(\mathscr{D}_{proj}'(\RR^{d}))$, the sequence $(m_n)_{n\in\NN_0}$ is also log-convex. 

We will show that there exists a finite positive constant $c_P$ such that 
\begin{equation}\label{boundCarl} 
\tilde{m}_n:= \sqrt{\sup_{f_1,\ldots,f_{2n}\in E}\left|\langle f_1\otimes\cdots\otimes f_{2n}, (_Pm)^{(2n)}\rangle\right|}\leq \sqrt{c_P m_{2n}}.
\end{equation}
The latter bound is sufficient to prove that the sequence $_Pm$ is determining. In fact, the log-convexity of $(m_n)_{n\in\NN_0}$ and the quasi-analyticity of $C\{{m}_n\}$ imply that the class $C\{\sqrt{c_P m_{2n}}\}$ is also quasi-analytic (see Lemma~\ref{CarlemanPari} and Proposition~\ref{ConstTimesCarl}). Hence,~\eqref{boundCarl} gives that $C\{\tilde{m}_n\}$ is also quasi-analytic.\\

It remains to show the bound in \eqref{boundCarl}.\\
Let us fix $n\in\NN$. Using Definition~\ref{SHIFT} and the assumption that $m$ is realized by $\mu$ on $\mathscr{D}'_{proj}(\RR^d)$, we get that for any $f_1,\ldots,f_{2n}\in\C_c^\infty(\RR^d)$ 
\begin{eqnarray*}
\left|\langle f_1\otimes\cdots\otimes f_{2n}, (_Pm)^{(2n)}\rangle\right|
\!\!&\leq&\!\!\sum_{j=0}^N\left|\int_{\mathscr{D}'_{proj}(\RR^d)}\langle p^{(j)}, \eta^{\otimes j}\rangle\langle f_1\otimes\cdots\otimes f_{2n},\eta^{\otimes (2n)} \rangle\mu(d\eta)\right|\\
\!\!&\leq&\!\!c_P\left(\int_{\mathscr{D}'_{proj}(\RR^d)}\left|\langle f_1\otimes\cdots\otimes f_{2n},\eta^{\otimes 2n} \rangle\right|^2\mu(d\eta)\right)^{\frac 12}\\
&=&c_P\left|\langle f_1^{\otimes 2}\otimes\cdots\otimes f_{2n}^{\otimes 2},m^{(4n)} \rangle\right|^{\frac 12},
\end{eqnarray*}
where 
$$c_P:=\sum_{j=0}^N\left(\int_{\mathscr{D}'_{proj}(\RR^d)}\left|\langle p^{(j)}, \eta^{\otimes j}\rangle\right|^2\mu(d\eta)\right)^{\frac 12}.$$
Note that $c_P$ is a finite positive constant since the realizing measure $\mu$ has finite local moments of any order. Hence, using the definition of $m_n$ and $\tilde{m}_n$, we get \eqref{boundCarl}.\\
%Consequently, we have that for any $f_1,\ldots,f_{2n}\in E$ 
%$$
%\left|\langle f_1\otimes\cdots\otimes f_{2n}, (_Pm)^{(2n)}\rangle\right|\leq c_P \sqrt{\left|\langle f_1^{\otimes 2}\otimes\cdots\otimes f_{2n}^{\otimes 2},m^{(4n)} \rangle\right|}\leq c_P m_{2n}.\\
%$$
\endproof
\vspace{-1cm}
\paragraph{\proof (Theorem~\ref{MainThm1})} \ \\ 
\textbf{Necessity}\ \\
Assume that $m$ is realized on $\mathcal{S}$ by a non-negative measure $\mu\in\mathcal{M}^*(\mathcal{S})$. Using \eqref{intFormFunct}, we get that for any $h\in\mathscr{P}_{\C_c^\infty}\left(\mathscr{D}_{proj}'(\RR^d)\right)$ and for any $i\in Y$ the following hold
$$L_m(h^2)=\int_{\mathcal{S}}{h^2(\eta)\,\mu(d\eta)}\quad\text{and}\quad L_m(P_ih^2)=\int_{\mathcal{S}}{P_i(\eta)h^2(\eta)\,\mu(d\eta)}.$$
Since integrals of non-negative functions w.r.t.\!\! a non-negative measure are non-negative, the inequalities in \eqref{MainThmCond1-1} hold. \\

\noindent\textbf{Sufficiency}\ \\
As already observed in Remark~\ref{semidefFunct}, the assumptions in \eqref{MainThmCond1-1} mean that the sequences $m$ and ${_P}m$ are positive semidefinite. Since $m$ is assumed to be determining, Theorem~\ref{KondrThm} guarantees the existence of a unique non-negative measure $\mu\in\mathcal{M}^*(\mathscr{D}_{proj}'(\RR^d))$ realizing $m$. On the one hand, according to Lemma~\ref{Prop2-1} the sequence $_{P_i} m$ is realized by the signed measure $P_i\mu$, i.e.\! for any $f^{(n)}\in\C_c^\infty(\RR^{nd})$ 
 \begin{equation}\label{eqA1}
 \langle f^{(n)}, (_{P_i}m)^{(n)}\rangle
  =\int_{\mathscr{D}_{proj}'(\RR^d)} \langle f^{(n)}, \eta^{\otimes n} \rangle P_i(\eta)\mu(d\eta).
 \end{equation}
On the other hand, by Proposition~\ref{Prop3-1}, the sequence $_{P_i}m$ is also determining. Hence, applying again Theorem~\ref{KondrThm}, the sequence $_{P_i}m$ is realized by a unique non-negative measure $\nu\in\mathcal{M}^*(\mathscr{D}_{proj}'(\RR^d))$, namely for any $f^{(n)}\in\C_c^\infty(\RR^{nd})$
\begin{equation}\label{eqB1}
 \langle f^{(n)}, (_{P_i}m)^{(n)}\rangle=\int_{\mathscr{D}_{proj}'(\RR^d)} \langle f^{(n)}, \eta^{\otimes n} \rangle \nu(d\eta).
 \end{equation}
Let $A_i:=\left\{\eta\in\mathscr{D}_{proj}'(\RR^d): P_i(\eta)\geq 0\right\}$ and let us define $\mu_i^+(B):=\mu(B\cap A_i)$ and $\mu_i^-(B):=\mu(B\cap (\mathscr{D}_{proj}'(\RR^d)\setminus A_i))$, for all $B\in\mathcal{B}(\mathscr{D}_{proj}'(\RR^d))$. Moreover, let us consider the non-negative measures $\sigma_i^+$ and $\sigma_i^-$ given by $\sigma_i^+(B):=\int_B P_i(\eta)\mu_i^+(d\eta)$ and $\sigma_i^-(B):=-\int_B P_i(\eta)\mu_i^-(d\eta)$, for all $B\in\mathcal{B}(\mathscr{D}_{proj}'(\RR^d))$.  Hence, we have that $\mu=\mu_i^{+}+\mu_i^-$ and $P_i\mu=\sigma_i^+-\sigma_i^-$. According to this notation, \eqref{eqA1} and~\eqref{eqB1} can be rewritten as
\begin{equation}\label{precedente}
\int\limits_{\mathscr{D}_{proj}'(\RR^d)} \!\!\langle f^{(n)}, \eta^{\otimes n} \rangle \sigma_i^+(d\eta)=\int\limits_{\mathscr{D}_{proj}'(\RR^d)}\!\!\langle f^{(n)}, \eta^{\otimes n} \rangle \sigma_i^-(d\eta)+\!\!\int\limits_{\mathscr{D}_{proj}'(\RR^d)}\!\!\langle f^{(n)}, \eta^{\otimes n} \rangle \nu(d\eta). 
\end{equation}
Since $m$ is determining and since $\mu^+\leq\mu$, the sequence $m^+$ consisting of all moment functions of $\mu^+$ is also determining. By Proposition~\ref{Prop3-1}, the sequence $_{P_i}{m^+}$ is determining, too.\\
As the two non-negative measures $\sigma_i^+$ and $\sigma_i^-+\nu$ both realize the determining sequence $_{P_i}{m^+}$, they coincide because Theorem~\ref{KondrThm} also guarantees the uniqueness of the realizing measure. 
This implies that the signed measure $P_i\mu$ is actually a non-negative measure on~$\mathscr{D}_{proj}'(\RR^d)$ and therefore, we have that 
\begin{equation}\label{Ai}
\forall\, i\in Y,\quad\mu\left(\mathscr{D}_{proj}'(\RR^d)\setminus A_i\right)=0.
\end{equation}
The set $\mathcal{S}=\bigcap_{i\in Y} A_i\in\sigma(\tau_w^{ind})$ by Corollary~\ref{S-Meas} and hence, $\mathcal{S}\in\sigma(\tau_w^{proj})$ by~\eqref{CorSigmaAlgebra}. It remains to show that $\mu$ is concentrated on the set~$\mathcal{S}$, i.e.\! $\mu\left(\mathscr{D}_{proj}'(\RR^d)\setminus\mathcal{S}\right)=0$. If $Y$ is countable, then the conclusion immediately follows from \eqref{Ai} using the countable subadditivity of $\mu$. In the case when $Y$ is uncountable, the latter argument does not work anymore but we can still get that the measure is concentrated on $\mathcal{S}$ proceeding as follows. First, let us extend $\mu$ to a measure $\mu'$ on $\mathscr{D}'_{ind}(\RR^d)$ by defining $\mu'(M):=\mu(M\cap\mathscr{D}'_{proj}(\RR^d))$, for all $M\in\sigma(\tau_w^{ind})$.
As $(\mathscr{D}'_{ind}(\RR^d), \tau_w^{ind})$ is a Radon space (see~Proposition~\ref{D'Radon}), the finite measure $\mu'$ is inner regular. This means that for any $M\in\sigma(\tau_w^{ind})$ and for any $\varepsilon>0$ there exists a compact set $K_\varepsilon\in\sigma(\tau_w^{ind})$ such that 
$
K_\varepsilon\subseteq M,
$
with
\begin{equation}\label{propertyB1}
\mu'(M)<\mu'( K_\varepsilon)+\varepsilon.
\end{equation}
Let us apply this property to $M= \mathscr{D}_{ind}'(\RR^d)\setminus\mathcal{S}=\bigcup_{i\in Y}\left(\mathscr{D}_{ind}'(\RR^d)\setminus A_i\right).$ 
Since the sets $\mathscr{D}_{ind}'(\RR^d)\setminus A_i$ form an open cover of $K_\varepsilon$, the compactness of $K_\varepsilon$ in $\left(\mathscr{D}_{ind}'(\RR^d), \tau^{ind}_w\right)$ implies that there exists a finite open subcover of $K_\varepsilon$, i.e. there exists a finite subset $J\subset Y$ such that
$K_\varepsilon\subseteq \bigcup_{i\in J} \left(\mathscr{D}_{ind}'(\RR^d)\setminus A_i\right).$
Therefore, we have that
$$
0\leq \mu'(K_\varepsilon)\leq \mu'\left(\bigcup_{i\in J} \left(\mathscr{D}_{ind}'(\RR^d)\setminus A_i\right)\right)
                                        \leq \sum_{i\in J} \mu\left(\left(\mathscr{D}_{ind}'(\RR^d)\setminus A_i\right)\cap\mathscr{D}_{proj}'(\RR^d)\right)=0,
$$
where in the last equality we used \eqref{Ai}.
Moreover, by  \eqref{propertyB1}, we have that $$\mu'\left(\mathscr{D}_{ind}'(\RR^d)\setminus\mathcal{S}\right)\leq\mu'(K_\varepsilon)+\varepsilon=\varepsilon.$$
Since this holds for any $\varepsilon>0$, we get $\mu'\left(\mathscr{D}_{ind}'(\RR^d)\setminus\mathcal{S}\right)=0$ and hence, we have $0=\mu'\left(\mathscr{D}_{ind}'(\RR^d)\setminus\mathcal{S}\right)=\mu\left((\mathscr{D}_{ind}'(\RR^d)\setminus\mathcal{S})\cap\mathscr{D}_{proj}'(\RR^d)\right)=\mu\left(\mathscr{D}_{proj}'(\RR^d)\setminus\mathcal{S}\right).$\\
\endproof

\begin{remark}\label{RemD'Ind}\ \\
Theorem~\ref{MainThm1} does still hold for any basic semi-algebraic set $\mathcal{S}$ which is subset of $\mathscr{D}_{ind}'(\RR^d)$ (instead of $\mathscr{D}_{proj}'(\RR^d)$) and gives a realizing measure actually concentrated on $\mathcal{S}\cap\mathscr{D}_{proj}'(\RR^d)$. If $\mathcal{S}\cap\mathscr{D}_{proj}'(\RR^d)=\emptyset$, then there is no contradiction because Theorem~\ref{MainThm1} shows that the only realizing measure is identically equal to zero, and so we know \emph{a posteriori} that all the moment functions were zeros. However, the case $\mathcal{S}\cap\mathscr{D}_{proj}'(\RR^d)\neq\emptyset$ is very common, since $\mathscr{D}_{proj}'(\RR^d)$ contains all tempered distributions, Radon measures and all locally integrable functions. Hence, if at least a single one of such generalized functions is contained in $\mathcal{S}$ then $\mathcal{S}\cap\mathscr{D}_{proj}'(\RR^d)\neq\emptyset$ and Theorem \ref{MainThm1} can be applied to get a non-zero realizing measure supported on $\mathcal{S}$, indeed on $\mathcal{S}\cap\mathscr{D}_{proj}'(\RR^d)$. Note that in Theorem~\ref{MainThm1} it is not sufficient to just assume that $m\in\mathcal{F}\left(\mathscr{D}'_{ind}(\RR^{d})\right)$. However, the assumption $m\in\mathcal{F}\left(\mathscr{D}'_{proj}(\RR^{d})\right)$ is not a restrictive requirement in any application. 
\end{remark}

\section{Applications}\label{Sec-Appl}
In this section we give some concrete applications of Theorem~\ref{MainThm1}.\\ 
In Subsection~\ref{FinDimCase}, we present Theorem~\ref{MainThm1} in the finite dimensional case. This theorem generalizes the results already know in literature about the classical moment problem on a basic semi-algebraic set of $\RR^d$.\\
In Subsection~\ref{1Subsec-Appl}, we study the case when we assume more regularity of type~\ref{condreg} on the putative moment functions, that is, we require that they are non-negative symmetric Radon measures. The advantage of this additional assumption is that it allows us to simplify the condition of determinacy and hence, to give an adapted version of Theorem~\ref{MainThm1}. 
In Subsection~\ref{ExRadonMeas}, we derive conditions on the putative moment functions to be realized by a random measure, that is, we assume $\mathcal{S}$ to be the set of all Radon measures on $\RR^d$. In this case, the fact that all the moment functions are themselves Radon measures is a necessary condition and so the results of Subsection~\ref{1Subsec-Appl} can be exploited. In Subsection~\ref{ExDensity}, we consider the case when $\mathcal{S}$ is the set of Radon measures with Radon-Nikodym densities w.r.t.\! the Lebesgue measure fulfilling an \emph{a priori} $L^\infty$ bound.

From now on let us denote by $\mathcal{R}(\RR^d)$ the space of all Radon measures on~$\RR^d$, namely the space of all non-negative Borel measures that are finite on compact sets in $\RR^d$.

\subsection{Finite dimensional case}\label{FinDimCase}\ \\
The $d-$dimensional moment problem on a closed basic semi-algebraic set $\mathcal{S}$ of $\RR^d$ is a special case of realizability problem. Hence, an analogous of Theorem~\ref{MainThm1} can be proved also in the finite dimensional case, where the condition $m:=(m^{(n)})_{n\in\NN_0}\in\mathcal{F}\left(\RR^{d}\right)$ holds for any multi-sequence of real numbers. In fact, if we denote by $\{e_1,\ldots, e_d\}$ the canonical basis of $\RR^d$ then we have that for each $n\in\NN_0$, $$m^{(n)}:=\sum\limits_{\stackrel{n_1,\ldots,n_d\in\NN_0}{n_1+\cdots+n_d=n}}m_{n_1,\ldots, n_d}^{(n)} \underbrace{e_1\otimes\cdots\otimes e_1}_{n_1\text{ times}}\otimes\cdots\otimes \underbrace{e_{d}\otimes\cdots\otimes e_d}_{n_d\text{ times}}\in\RR^{dn}.$$ 
%$m^{(n)}:=\sum\limits_{j_1,\ldots,j_n=1}^d m_{j_1,\ldots, j_n}^{(n)} e_{j_1}\otimes\cdots\otimes e_{j_n}\in\RR^{dn}$. 
The notion of polynomials, quadratic module and Riesz's functional given at the beginning of Section~\ref{Sec-Main}, in the $d-$dimensional case coincide with the classical ones. \\
The condition of determinacy on $m$ reduces to the requirement that the class $C\left\{\sqrt{\max\limits_{\stackrel{n_1,\ldots,n_d\in\NN_0}{n_1+\cdots+n_d=2n}}|m_{n_1,\ldots, n_d}^{(2n)}|}\right\}$ is quasi-analytic. This follows by taking the subset $E:=\{e_1,\ldots,e_d\}$ in Definition~\ref{DefSeq}.\\ 
In this framework, the whole proof we made in the infinite dimensional case can be employed as well, taking in consideration that $\RR^d$ is Polish and so Radon. Actually, we can even get a stronger result by refining our proof in finite dimensions. Indeed, if we replace the assumption of $m$ being determining with the classical multivariate Carleman condition, that is for any $i\in\{1,\ldots,d\}$ the class $C\left\{\sqrt{|m_{0,\ldots,0, 2n,0,\ldots, 0}^{(2n)}|}\right\}$ is quasi-analytic (where $2n$ is at the $i-$th position of the index $d-$tuple), then we can still use the same proof but we need to substitute Theorem~\ref{KondrThm} with the $d-$dimensional version of Hamburger's theorem (see e.g.\! \cite{Sh-Tam43, Nuss65, Berg87}). In this way, we obtain the following general result.
 \begin{theorem}\label{MainThm-FinDimCase}\ \\
Let $m$ be a multi-sequence of real numbers, which fulfills the classical multivariate Carleman condition and let $$
\mathcal{S}=\bigcap_{i\in Y}\left\{\r\in\RR^d| \ P_i(\r)\geq 0\right\},
$$
where $Y$ is an index set not necessarily countable and $P_i\in\mathscr{P}_{\RR}\left(\RR^d\right)$ that is polynomial on $\RR^d$ with real coefficients. Then $m$ is realized by a unique non-negative measure $\mu\in\mathcal{M}^*(\mathcal{S})$ if and only if the following inequalities hold
$$
L_m(h^2)\geq 0, \,\,L_m(P_i h^2)\geq 0\,\,,\,\, \forall h\in\mathscr{P}_{\RR}\left(\RR^d\right),\,\forall i\in Y.
$$
Equivalently, if and only if the functional $L_m$ is non-negative on the quadratic module $\mathcal{Q}(\mathscr{P}_{\mathcal{S}})$.
\end{theorem}

This theorem extends the result given by Lasserre in~\cite{Las2011}. In fact, Theorem~\ref{MainThm-FinDimCase} includes the case when $\mathcal{S}$ is defined by an uncountable family of polynomials. Furthermore, the classical multivariate Carleman condition assumed in Theorem~\ref{MainThm-FinDimCase} is a more general bound than the one assumed in \cite{Las2011}.

\subsection{Realizability of Radon measures}\label{1Subsec-Appl}
\begin{definition}\label{WeighCarl}\ \\
A sequence $m\in\mathcal{F}\left(\mathcal{R}(\RR^{d})\right)$ satisfies the \emph{weighted Carleman type condition} if for each $n\in\NN$ there exists a function $k^{(n)}_2\in\C^\infty(\RR^d)$ with $k^{(n)}_2(\r)\geq 1$ for all $\r\in\RR^d$ such that
\begin{equation}\label{CarlGen}
 \sum_{n=1}^{\infty}\frac{1}{\left(\sup\limits_{\stackrel{\z\in\RR^d}{\|\z\|\leq n}}\sup\limits_{\x \in [-1,1]^d}\sqrt{ \tilde{k}_2^{(n)}(\z+\x)}\right)\!\!\sqrt[2n]{\int_{\RR^{2nd}}\frac{m^{(2n)}(d\r_1,\ldots,d\r_{2n})}{\prod_{l=1}^{2n}k^{(2n)}_2(\r_l)}}}=\infty,
 \end{equation}  
where $\tilde{k}^{(n)}_2\in\C^\infty(\RR^d)$ such that $\tilde{k}^{(n)}_2(\r)\geq\left|(D^\kappa k^{(n)}_2)(\r)\right|^2$ for all $|\kappa|\leq \lceil{\frac {d+1}{2}}\rceil$.
\end{definition}
As suggested by the name, the condition \eqref{CarlGen} is an infinite dimensional weighted version of the classical Carleman condition, which ensures the uniqueness of the solution to the $d-$dimensional moment problem (for $d=1$ see \cite{Carl26}, for $d\geq 2$ see e.g.~\cite{Sh-Tam43, Nuss65, Berg87, DeJeu03}) .

\begin{corollary}\label{MainThm2}\ \\
Let $m\in\mathcal{F}\left(\mathcal{R}(\RR^{d})\right)$ fulfill the weighted Carleman type condition in Definition~\ref{WeighCarl} and let $\mathcal{S}\subseteq\mathscr{D}'(\RR^d)$ be a basic semi-algebraic set of the form~\eqref{S-Semialg}.
Then $m$ is realized by a unique non-negative measure $\mu\in\mathcal{M}^*(\mathcal{S})$ with
\begin{equation}\label{FiniteGlobalMom}
\int_{\mathcal{S}}\langle \frac{1}{k^{(n)}_2},\eta\rangle^n \mu(d\eta)<\infty,\quad \forall\,n\in\NN_0,
\end{equation}
if and only if the following inequalities hold
\begin{equation}\label{MainThmCond1}
L_m(h^2)\geq 0, \,\,L_m(P_i h^2)\geq 0,\,\,\, \forall h\in\mathscr{P}_{\C_c^\infty}\left(\mathscr{D}_{proj}'(\RR^d)\right),\,\forall i\in Y,
\end{equation}
and for any $n\in\NN_0$ we have
\begin{equation}\label{MainThmCond2}
\int_{\RR^{2nd}}\frac{m^{(2n)}(d\r_1,\ldots,d\r_{2n})}{\prod_{l=1}^{2n}k^{(2n)}_2(\r_l)}<\infty.
\end{equation}
\end{corollary}

\begin{remark}\label{LemmaFinite1}\ \\
If $m$ is realized by a non-negative measure $\mu\in\mathcal{M}^*(\mathscr{D}_{proj}'(\RR^d))$ and $m$ satisfies \eqref{CarlGen} then \eqref{MainThmCond2} holds also for the odd orders.
\end{remark}
Corollary~\ref{MainThm2} is essentially a consequence of the following proposition.
 \begin{prop}\label{Prop1-1}\ \\
If $m$ satisfies \eqref{CarlGen} and \eqref{MainThmCond2}, then $m$ is a determining sequence in the sense of Definition~\ref{DefSeq}.
 \end{prop}
\proof\ \\
Let us preliminarily recall that $\mathcal{R}(\RR^d)\subset\mathscr{D}'(\RR^d)$ and so $m$ is automatically in $\mathcal{F}(\mathscr{D}'(\RR^d))$ as required by Definition~\ref{DefSeq}.\\
For any $f_1,\ldots,f_{n}\in\C_c^\infty(\RR^d)$ and any $n\in\NN$ we can easily see that
\begin{equation}\label{PrevStima}
 \left|\left\langle f_1\otimes\cdots\otimes f_{n}, m^{(n)} \right\rangle\right|\leq\int_{\RR^{nd}}\prod_{l=1}^{n}k^{(n)}_2(\r_l)\left|f_l(\r_l)\right|\frac{m^{(n)}(d\r_1,\ldots,d\r_{n})}{\prod_{l=1}^{n}k^{(n)}_2(\r_l)}.
\end{equation}
By the Sobolev embedding theorem for weighted spaces (see~\cite{B86}), we get that for any $\tilde{k}^{(n)}_2\in\C^\infty(\RR^d)$ with $\tilde{k}^{(n)}_2(\r)\geq\left|(D^\kappa k^{(n)}_2)(\r)\right|^2$ for all $|\kappa|\leq \big\lceil{\frac {d+1}{2}}\big\rceil$, $\C_c(\RR^d)\subseteq H_{ \tilde{k}^{(n)}}$, where ${H_{ \tilde{k}^{(n)}}}:={W_2^{\lceil{\frac {d+1}{2}}\rceil}(\RR^d, \tilde{k}^{(n)}_2(\r)d\r)}$ and $\tilde{k}^{(n)}:=\left(\big\lceil{\frac {d+1}{2}}\big\rceil, \tilde{k}^{(n)}_2\right)$. Using this result in \eqref{PrevStima}, we have that there exists a finite positive constant $C$ such that 
$$ \left|\left\langle f_1\otimes\cdots\otimes f_{n}, m^{(n)} \right\rangle\right|\leq\ C^{n}\prod_{l=1}^{n}\left\|f_l(\r_l)\right\|_{H_{ \tilde{k}^{(n)}}}\int\limits_{\RR^{nd}}\frac{m^{(n)}(d\r_1,\ldots,d\r_{n})}{\prod_{l=1}^{n}k^{(n)}_2(\r_l)}.
 $$
Hence, by choosing $E$ as in Lemma~\ref{LemmaE}, we have that 
\begin{eqnarray}
m_n\!\!\!\!\!&:=&\!\!\!\!\!\sqrt{\sup_{f_1,\ldots,f_{2n}\in E}\left|\langle f_1\otimes\cdots\otimes f_{2n}, m^{(2n)}\rangle\right|}\nonumber\\
&\leq&\!\!\!\!\!\sqrt{C^{2n}\left(\sup_{f\in E}\left\|f\right\|_{H_{ \tilde{k}^{(n)}}}\right)^{2n}\int\limits_{\RR^{2nd}}\frac{m^{(2n)}(d\r_1,\ldots,d\r_{n})}{\prod_{l=1}^{2n}k^{(2n)}_2(\r_l)}}\nonumber\\
&\leq&\!\!\!\!\!\!\left(Cc^d_{\lceil{\frac {d+1}{2}}\rceil}\sup\limits_{\stackrel{\z\in\RR^d}{\|\r\|\leq n}}\!\!\sup\limits_{\x \in [-1,1]^d}\sqrt{ \tilde{k}_2^{(2n)}(\z+\x)}\right)^n\!\!\!\sqrt{\int\limits_{\RR^{2nd}}\!\frac{m^{(2n)}(d\r_1,\ldots,d\r_{2n})}{\prod_{l=1}^{2n}k^{(2n)}_2(\r_l)}}\label{Bound-mn}.
\end{eqnarray}
Then the condition \eqref{MainThmCond2} guarantees that the $m_n$'s are finite and~\eqref{CarlGen} implies that the class $C\{m_n\}$ is quasi-analytic.\\
\endproof
\paragraph{\proof (Corollary~\ref{MainThm2})} \ \\
Since the necessity part follows straightforwardly, let us focus on the sufficiency.\\
%\textbf{Necessity}\ \\
%Assume that $m$ is realized by a non-negative measure $\mu\in\mathcal{M}^*(\mathcal{S})$. Then, by Theorem~\ref{MainThm1}, the inequalities \eqref{MainThmCond1} are fulfilled and \eqref{MainThmCond2} follows by Remark~\ref{LemmaFinite1}.\\
%\textbf{Sufficiency}\ \\
Since $m$ is determining by Proposition \ref{Prop1-1} and \eqref{MainThmCond1} holds by assumption, we can apply Theorem~\ref{MainThm1} to get that $m$ is realized by $\mu\in\mathcal{M}^*(\mathcal{S})$.

It remains to show \eqref{FiniteGlobalMom}. For any positive real number $R$ let us define a function $\chi_R$ such that
\begin{equation}\label{ChiR}
\chi_R\in\C_c^\infty(\RR^d)\, \text{ and }\, \chi_R(\r):=\left\{\begin{array}{ll}
1 & \text{if } |\r|\leq R\\
0 & \text{if } |\r|\geq R+1.
\end{array}
\right. 
\end{equation}
Since $m$ is realized by $\mu\in\mathcal{M}^*(\mathcal{S})$, for any $n\in\NN_0$ and for any positive real number $R$ we have that
$$\int_{\mathcal{S}}\langle\frac{\chi_R}{k^{(n)}_2},\eta\rangle^{ n}\mu(d\eta)=\int_{\RR^{nd}}\prod_{l=1}^{n}\frac{\chi_R(\r_l)}{k^{(n)}_2(\r_l)}m^{(n)}(d\r_1,\ldots,d\r_{n}).$$
Hence, the monotone convergence theorem for $R\to\infty$ and Remark~\ref{LemmaFinite1} give \eqref{FiniteGlobalMom}.\\
\endproof

\begin{remark}\label{RemarkNorm}\ \\
The proof of Proposition \ref{Prop1-1} is a particular instance of what we were pointing out in Remark \ref{RemarkE}. In fact, the regularity assumed on the sequence $m$, that is $m$ consisting of Radon measures, allowed us to get the bound~\eqref{Bound-mn} from \eqref{CarlGen} and \eqref{MainThmCond2} for some index $\tilde{k}^{(n)}=(\tilde{k}_1^{(n)},\tilde{k}_2^{(n)})$ with $\tilde{k}_1^{(n)}=\big\lceil\frac{d+1}{2}\big\rceil$ and so independent of $n$. \\ 
Note that to obtain this result it was important to use our definition of determining sequence (see Definition~\ref{DefSeq}). In fact, if we used the one given in \cite{BeKo88} involving the norms $\| m^{(2n)}\|_{H_{-k^{(2n)}}^{\otimes 2n}}$ (see Remark~\ref{RemarkYuri}), we would have got $\tilde{k}_1^{(n)}>\big\lceil\frac{n(d+1)}{2}\big\rceil$ and as a consequence an extra factor of at least order $(2n)!$ under the root in~\eqref{CarlGen}. This observation is in line with Remark 3 in \cite[Vol.~II, p.73]{BeKo88}.
\end{remark}

If we assume even more regularity on $m$, then Corollary~\ref{MainThm2} takes the following simpler form.
\begin{corollary}\label{MainThm2-simpler}\ \\
Let $m\in\mathcal{F}\left(\mathcal{R}(\RR^{d})\right)$ be such that for some $k_2\in\C^\infty(\RR^d)$, independent of $n$, with $k_2(\r)\geq 1$ for all $\r\in\RR^d$ the following holds
\begin{equation*}\label{CarlGen1}
 \sum_{n=1}^{\infty}\frac{1}{\sqrt[2n]{\int_{\RR^{2nd}}\frac{m^{(2n)}(d\r_1,\ldots,d\r_{2n})}{\prod_{l=1}^{2n}k_2(\r_l)}}}=\infty.
 \end{equation*} If $\mathcal{S}\subseteq\mathscr{D}'(\RR^d)$ is a basic semi-algebraic set of the form~\eqref{S-Semialg}, then $m$ is realized by a unique non-negative measure $\mu\in\mathcal{M}^*(\mathcal{S})$ with
$$
\int_{\mathcal{S}}\langle \frac{1}{k_2},\eta\rangle^n \mu(d\eta)<\infty,\quad \forall\,n\in\NN_0,
$$
if and only if the following inequalities hold
$$
L_m(h^2)\geq 0, \,\,L_m(P_i h^2)\geq 0,\,\,\, \forall h\in\mathscr{P}_{\C_c^\infty}\left(\mathscr{D}_{proj}'(\RR^d)\right),\,\forall i\in Y,
$$
and for any $n\in\NN_0$ we have
$$
\int_{\RR^{2nd}}\frac{m^{(2n)}(d\r_1,\ldots,d\r_{2n})}{\prod_{l=1}^{2n}k_2(\r_l)}<\infty.
$$
\end{corollary}

\subsection{Realizability on the space of Radon measures $\mathcal{R}(\RR^d)$}\label{ExRadonMeas}
\begin{example}\label{ApplRadon}\ \\
The set $\mathcal{R}(\RR^d)$ of all Radon measures on $\RR^d$ is a basic semi-algebraic subset of $\mathscr{D}'(\RR^d)$, i.e.
\begin{equation}\label{radonMeasRepr}
\mathcal{R}(\RR^d)=\bigcap_{\varphi\in\C_c^{+,\infty}(\RR^d)}\left\{\eta\in\mathscr{D}'(\RR^d): \Phi_\varphi(\eta)\geq 0\right\}
\end{equation}
where $\Phi_\varphi(\eta):=\langle \varphi, \eta \rangle$. 
\end{example}
\proof\ \\
The representation \eqref{radonMeasRepr} follows from the fact that there exists a one-to-one correspondence between the Radon measures on $\RR^d$ and the continuous non-negative linear functionals on the space $\mathscr{D}_{proj}(\RR^d)$. In fact, for any $\eta\in\mathcal{R}(\RR^d)$ the functional 
\begin{eqnarray*}
\C_c^\infty(\RR^d)&\to&\RR \\
 \varphi & \mapsto&\langle \varphi,\eta\rangle=\int_{\RR^d} \varphi(\r)\eta( d\r)
 \end{eqnarray*}
is non-negative and it is an element of $\mathscr{D}'(\RR^d)$. Conversely, by a theorem due to L.~Schwartz (see \cite[Theorem V]{Sch57} %\cite[Theorem 5.3.1]{BlaBru03}
), every non-negative linear functional on $\C_c^\infty(\RR^d)$ can be represented as integral w.r.t.\! a Radon measure on $\RR^d$. \\
\endproof

Using the representation \eqref{radonMeasRepr}, we obtain a realizability theorem for $\mathcal{S}=\mathcal{R}(\RR^d)$, namely  Corollary~\ref{MainThm2} becomes
\begin{theorem}\label{ThmRadonMeasures}\ \\
Let $m\in\mathcal{F}\left(\mathcal{R}(\RR^{d})\right)$ fulfill the weighted Carleman type condition \eqref{CarlGen}. Then $m$ is realized by a unique non-negative measure $\mu\in\mathcal{M}^*(\mathcal{R}(\RR^d))$ with
$$
\int_{\mathcal{S}}\langle \frac{1}{k_2^{(n)}},\eta\rangle^n \mu(d\eta)<\infty,\quad \forall\,n\in\NN_0,
$$
if and only if the following inequalities hold
\begin{eqnarray}
&\ & L_m(h^2)\geq 0\,,\,\, \forall h\in\mathscr{P}_{\C_c^\infty}\left(\mathscr{D}_{proj}'(\RR^d)\right),\label{Semid}\\
&\ & L_m(\Phi_\varphi h^2)\geq 0\,,\,\, \forall h\in\mathscr{P}_{\C_c^\infty}\left(\mathscr{D}_{proj}'(\RR^d)\right),\, \forall \varphi\in\C_c^{+,\infty}(\RR^d),\label{Semid1}\\
&\ & \int_{\RR^{2nd}}\frac{m^{(2n)}(d\r_1,\ldots,d\r_{2n})}{\prod_{l=1}^{2n}k_2^{(2n)}(\r_l)}<\infty,\,\,\forall n\in\NN_0.\label{qxrt}\end{eqnarray}
\end{theorem}
Note that if $\mu$ is concentrated on $\mathcal{R}(\RR^d)$ then $m_\mu^{(n)}\in\mathcal{R}(\RR^{dn})$ for all $n\in\NN_0$.\\

The previous theorem still holds even when $m$ does not consist of Radon measures. In this case, instead of \eqref{CarlGen} and \eqref{qxrt}, one has to assume that $m$ is determining in the sense of Definition \ref{DefSeq}

The assumption \eqref{CarlGen} can be actually weakened by taking into account a result due to S.N. \v{S}ifrin about the infinite dimensional moment problem on dual cones in nuclear spaces (see \cite{S74}). Indeed, applying \v{S}ifrin's results to the cone $\C_c^{+,\infty}(\RR^d)$, it is possible to obtain a particular instance of our Theorem~\ref{MainThm1} for the case $\mathcal{S}=\mathcal{R}(\RR^d)$ (the latter is in fact the dual cone of $\C_c^{+,\infty}(\RR^d)$) but with the difference that in the determinacy condition the quasi-analyticity of the $m_n$'s is replaced by the so-called Stieltjes condition $\sum_{n=1}^\infty  m_n^{-\frac{1}{2n}}=\infty$. As a consequence, the condition~\eqref{CarlGen} in Theorem~\ref{ThmRadonMeasures} can be replaced by the following weaker one
\begin{equation*}\label{StieltjesGen}
\sum_{n=1}^{\infty}\frac{1}{\sqrt{\sup\limits_{\stackrel{\z\in\RR^d}{\|\z\|\leq n}}\sup\limits_{\x \in [-1,1]^d}\sqrt{ \tilde{k}_2^{(n)}(\z+\x)}}\sqrt[4n]{\int_{\RR^{2nd}}\frac{m^{(2n)}(d\r_1,\ldots,d\r_{2n})}{\prod_{l=1}^{2n}k^{(2n)}_2(\r_l)}}}=\infty,
 \end{equation*}  
which we call \emph{weighted generalized Stieltjes condition}.
\begin{remark}\label{RemarkSDefCondMomMeas}\ \\
The condition \eqref{Semid} can be rewritten as
\footnotesize
$$\sum_{i,j}\langle h^{(i)}\otimes h^{(j)},\, m^{(i+j)}\rangle\geq 0,\quad\forall\, h^{(i)}\in\C_c^\infty(\RR^{id}),$$
\normalsize
and \eqref{Semid1} as
\footnotesize
$$\sum_{i,j}\langle h^{(i)}\otimes h^{(j)}\otimes\varphi,\, m^{(i+j+1)}\rangle\geq 0,\quad\forall\, h^{(i)}\in\C_c^\infty(\RR^{id}),\, \forall \varphi\in\C_c^{+,\infty}(\RR^d).$$
\normalsize
Recalling Definition~\ref{SHIFT}, we can restate these conditions as follows: the sequence $(m^{(n)})_{n\in\NN_0}$ and its shifted version $((_{\Phi_\varphi}m)^{(n)})_{n\in\NN_0}$ are positive semidefinite in the sense of Definition~\ref{PosSemiDef}.\\ In particular, if for each $n\in\NN_0$, $m^{(n)}$ has a Radon-Nikodym density, that is there exists $\alpha^{(n)}\in L^1(\RR^n,\lambda)$ s.t.\! $m^{(n)}(d\r_1,\ldots, d\r_n)=\alpha^{(n)}(\r_1,\ldots, \r_n)d\r_1\cdots d\r_n$, then \eqref{Semid} and \eqref{Semid1} can be rewritten as
$$
\resizebox{.9\hsize}{!}{$\sum\limits_{i,j}\int_{\RR^{d(i+j)}} h^{(i)}(\r_1,\ldots,\r_i) h^{(j)}(\r_{i+1},\ldots,\r_{i+j}) \alpha^{(i+j)}(\r_1,\ldots, \r_{i+j})d\r_1\cdots d\r_{i+j}\geq 0,$}$$
$$\resizebox{\hsize}{!}{$\sum\limits_{i,j}\int_{\RR^{d(i+j+1)}}h^{(i)}(\r_1,\ldots,\r_i) h^{(j)}(\r_{i+1},\ldots,\r_{i+j})\varphi(\y) \alpha^{(i+j+1)}(\r_1,\ldots, \r_{i+j}, \y)d\r_1\cdots d\r_{i+j}d \y\geq 0.$}$$

These conditions can be interpreted as that $(\alpha^{(n)})_{n\in\NN_0}$ is positive semidefinite and
that for $\lambda-$almost all $\y\in\RR^d$ the sequence $(\alpha^{(n+1)}(\cdot,\y))_{n\in\NN_0}$ is positive semidefinite, where the positive semidefiniteness is intended in a generalized sense. In this reformulation the analogy with the Stieltjes moment problem is evident, since necessary and sufficient conditions for the realizability on $\RR^+$ of a sequence of numbers $(m_n)_{n\in\NN_0}$ are that $(m_n)_{n\in\NN_0}$ and $(m_{n+1})_{n\in\NN_0}$ are positive semidefinite.
\end{remark}

The measure constructed in Theorem \ref{ThmRadonMeasures} lives on the Borel $\sigma-$algebra generated by the weak topology $\tau_w^{proj}$ on $\mathscr{D}'_{proj}$ restricted to its subset $\mathcal{R}(\RR^d)$. A natural topology on $\mathcal{R}(\RR^d)$ is the vague topology $\tau_v$, i.e.\! the smallest topology such that the mappings
$$
\eta \mapsto \langle f, \eta\rangle=\int_{\RR^d} f(\r)
\eta(d\r)
$$
are continuous for all $f \in \C_c(\RR^d)$. These two topologies actually coincide on $\mathcal{R}(\RR^d)$.

This result directly follows from the Hausdorff criterion %(see \cite[Theorem 4.8, p.35]{Willard}) 
if one intersects the neighbourhood bases with sets of the following form
\begin{equation*}
U_{\chi_\varphi; N}:=\left\{\eta\in\mathcal{R}(\RR^d) :  \left|\langle \chi_\varphi, \eta-\nu\rangle\right|<N \right\},
\end{equation*}
where $N$ is a positive integer and $\chi_\varphi$ is a smooth characteristic function of the support of a function $\varphi\in\C_c(\RR^d)$ (see \eqref{ChiR}).\\ As a consequence of the equivalence of the two topologies, the associated Borel $\sigma-$algebras also coincide and they are equal to $\sigma(\tau_w^{proj})\cap\mathcal{R}(\RR^d)$. 

\subsection{Realizability on the set of measures with bounded density}\label{ExDensity}
\begin{example}\ \\
Let $c\in\RR^+$. The set $\mathcal{S}_c$ of all Radon measures with density w.r.t. the Lebesgue measure $\lambda$ on $\RR^d$ which is $L^\infty-$bounded by $c$, i.e.
\begin{equation}\label{defACMeas}
\mathcal{S}_c:=\left\{\eta\in\mathcal{R}(\RR^d): \eta(d \r)=f(\r)\lambda(d\r)\,\,\text{with}\,\,f\geq 0\,\,\text{and}\,\,\|f\|_{L^\infty}\leq c\right\} 
\end{equation}
is a basic semi-algebraic subset of $\mathscr{D}'(\RR^d)$. More precisely, we get that
\begin{equation}\label{ACMeasSemiAlg}
\mathcal{S}_c=\mathcal{R}(\RR^d) \cap \bigcap_{\varphi\in\C_c^{+,\infty}(\RR^d)}\left\{\eta\in\mathscr{D}'(\RR^d): c\langle \varphi, \lambda\rangle-\langle \varphi, \eta\rangle\geq 0\right\}.
\end{equation} 
\end{example}
\proof\ \\
{\bf Step I}: $\subseteq$\ \\
Let $\eta\in\mathcal{S}_c$, then by definition \eqref{defACMeas}, we get that for any $\varphi\in\C_c^{+,\infty}(\RR^d)$ 
$$
\langle \varphi, \eta\rangle=\int_{\RR^d}\varphi(\r)f(\r)\lambda(d\r)\leq\|f\|_{L^\infty}\int_{\RR^d}\varphi(\r)\lambda(d\r)\leq c\langle \varphi, \lambda\rangle. 
$$

\noindent{\bf Step II}: $\supseteq$\ \\
Let $\eta\in\mathcal{R}(\RR^d)$ such that  
\begin{equation}\label{relazione}
c\langle \varphi, \lambda\rangle-\langle \varphi, \eta\rangle\geq 0,\,\,\forall\,\varphi\in\C_c^{+,\infty}(\RR^d).
\end{equation}
By density, the previous condition holds for all $\varphi\in L^1(\RR^d, \lambda-\eta)$ 
and in particular for $\varphi=\Ii_A$, where $A\in\mathcal{B}(\RR^d)$ bounded. Hence, $\eta\ll\lambda$ and so, by the Radon-Nikodym theorem, there exists $f\geq 0$ such that 
\begin{equation}\label{relazione2}
\eta(d \r)=f(\r)\lambda(d\r).
\end{equation}
By \eqref{relazione2} and \eqref{relazione}, for any $A\in\mathcal{B}(\RR^d)$ bounded we get that%$(\lambda-\eta)-$measurable set 
$$\int_Af(\r)\lambda(d\r)=\int_A\eta(d\r)\leq c\int_A\lambda(d\r).$$
Hence, $f(\r)\leq c$ $\lambda-$a.e. in each bounded $A$ and therefore $\|f\|_{L^\infty}\leq c$.\\
% $f(\r)\leq c$ $\lambda-$a.e. in each bounded $A$
%implies that in every $A$ bounded  there exists a set of measure zero such that $f(\r)> c$
%since \RR^d can be written as union of balls (which are bounded), then there are infinitely many of these null sets in \RR^d where $f(\r)> c$ and this means that the set given by the union of this null sets is still a null set such that $f(\r)> c$, i.e. 
% $f(\r)\leq c$ $\lambda-$a.e. in $\RR^d$
%Recall that $\|f\|_{L^\infty}:=\inf\{M\in\RR^+: |f(\r)|\leq M,\,\,\lambda-\text{a.e.}\}$  implies that $\|f\|_{L^\infty}\leq c$
\endproof

Using the representation \eqref{ACMeasSemiAlg}, we can explicitly rewrite Corollary~\ref{MainThm2} for $\mathcal{S}=\mathcal{S}_c$ as follows.
\begin{theorem}\label{ThmACMeas}\ \\
Let $c\in\RR^+$. Let $m\in\mathcal{F}\left(\mathcal{R}(\RR^{d})\right)$ fulfill the weighted Carleman type condition~\eqref{CarlGen}.
Then $m$ is realized by a unique non-negative measure $\mu\in\mathcal{M}^*(\mathcal{S}_c)$ with
$$
\int_{\mathcal{S}}\langle \frac{1}{k_2^{(n)}},\eta\rangle^n \mu(d\eta)<\infty,\quad \forall\,n\in\NN_0,
$$if and only if the following inequalities hold.
\begin{eqnarray}
&\ & L_m(h^2)\geq 0,\,\, \forall h\in\mathscr{P}_{\C_c^\infty}\left(\mathscr{D}_{proj}'(\RR^d)\right),\label{condizione1}\\
&\ & L_m(\Phi_\varphi h^2)\geq 0,\,\,\forall h\in\mathscr{P}_{\C_c^\infty}\left(\mathscr{D}_{proj}'(\RR^d)\right), \forall \varphi\in\C_c^{+,\infty}(\RR^d),\label{condizione2}\\
&\ & L_m(\Gamma_{c,\varphi} h^2)\geq 0,\,\,\forall h\in\mathscr{P}_{\C_c^\infty}\left(\mathscr{D}_{proj}'(\RR^d)\right), \forall \varphi\in\C_c^{+,\infty}(\RR^d),\label{condizione3}\\
& \ & \int_{\RR^{2nd}}\frac{m^{(2n)}(d\r_1,\ldots,d\r_{2n})}{\prod_{l=1}^{2n}k_2^{(2n)}(\r_l)}<\infty\,,\,\,\forall n\in\NN_0,\nonumber\label{condizione4}
\end{eqnarray}
where $\Phi_\varphi(\eta):=\langle \varphi, \eta \rangle$ and $\Gamma_{c,\varphi}(\eta):=c\langle \varphi,\lambda\rangle-\langle \varphi, \eta\rangle.$ 
\end{theorem}
\begin{remark}\label{RemarkSDefCondMomMeas2}\ \\
Proceeding as in Remark~\ref{RemarkSDefCondMomMeas}, we can work out the analogy between the realizability problem on $\mathcal{S}_c$ and the moment problem on $[0,c]$. Indeed, if each $m^{(n)}$ has density $\alpha^{(n)}$ w.r.t. the Lebesgue measure, then \eqref{condizione1}, \eqref{condizione2} and \eqref{condizione3} mean just that $(\alpha^{(n)})_{n\in\NN_0}$ is positive semidefinite and that, for $\lambda-$almost all $\y\in\RR^d$, $(\alpha^{(n+1)}(\cdot,\y))_{n\in\NN_0}$ and $(c\alpha^{(n)}(\cdot)-\alpha^{(n+1)}(\cdot,\y))_{n\in\NN_0}$ are positive semidefinite. 
Similarly, necessary and sufficient conditions for the realizability on $[0,c]$ of a sequence of numbers~$(m_n)_{n\in\NN_0}$, where
$$[0,c]=\{x\in\RR: x\geq 0\}\cap \{x\in\RR: c-x\geq 0\},$$
are that $(m_n)_{n\in\NN_0}$, $(m_{n+1})_{n\in\NN_0}$ and $(c\cdot m_n-m_{n+1})_{n\in\NN_0}$ are positive semidefinite (see \cite{Devi53} and \cite{BeMa}).
\end{remark}

\section{Appendix}\label{Sec-App}
\subsection{Quasi-analyticity}\label{Sec-App1}\ \\
Let us recall the basic definitions and state the results used throughout this paper concerning the theory of quasi-analyticity.
\begin{definition}[The class $C\{M_n\}$]\label{quasiAn}\ \\
Given a sequence of positive real numbers $(M_n)_{n\in\NN_0}$, we define the class $C\{M_n\}$ as the set of all functions $f\in\C^\infty(\RR)$ such that for any $n\in\NN_0$
$$\left\|D^n f\right\|_\infty\leq\beta_f B_f^n M_n,$$
where $D^n f$ is the $n-$th derivative of $f$, $\left\|D^n f\right\|_\infty:=\sup\limits_{x\in\RR}\left|D^n f(x)\right|$, and $\beta_f$, $B_f$ are positive constants only depending on $f$.
\end{definition}
\begin{definition}[Quasi-analytical class]\label{QuasianalityCLA}\ \\
A class $C\{M_n\}$ is said to be quasi-analytic if the conditions
$$f\in C\{M_n\},\,\,(D^n f)(0)=0,\quad\forall\, n\in\NN_0,$$
imply that $f(x)=0$ for all $x\in\RR$.
\end{definition}

The main result in the theory of quasi-analyticity is the Denjoy-Carleman theorem, which is easy to prove when the sequence is log-convex and has the first term equal to 1 (see~\cite{Ru74} for a proof of the theorem in this case). 
\begin{definition}[Log-convexity]\label{LogConv1}\ \\
A sequence of positive real numbers $(M_n)_{n\in\NN_0}$ is said to be log-convex if and only if for all $n \geq 1 $ we have that
$M_n^2 \leq M_{n-1} M_{n+1}$.
\end{definition}

However, when we deal with classes of functions, the assumption of log-convexity and the assumption $M_0=1$ actually involve no loss of generality. In fact, one can prove that for any sequence $(M_n)_{n\in\NN_0}$ there always exists a log-convex sequence $({M}^c_n)_{n\in\NN_0}$ such that the classes $C\{M_n\}$ and $C\{{M}^c_n\}$ coincide. More precisely, the sequence $({M}^c_n)_{n\in\NN_0}$ is the convex regularization of $({M}_n)_{n\in\NN_0}$ by means of the logarithm (for more details on this regularization see \cite{Mand52}).  % \cite[Chapter VI, Theorem 6.5.III]{Mand52}  il teorema sulle quasi analitiche and \cite{Gorny39} Gorny does it for trigonometric classes). \\
Hence, we have that $C\{M_n\}$ is quasi-analytic if and only if $C\{M^c_n\}$ is quasi-analytic (see \cite[Chapter VI, Theorem~6.5.III]{Mand52}). Clearly, if $(M_n)_{n\in\NN_0}$ is log-convex then $M_n^c\equiv M_n$ for all $n\in\NN_0$. Furthermore, if $M_0\neq 1$ then one can always normalize the sequence and consider $(\frac{M_n}{M_0})_{n\in\NN_0}$, since it is easy to see that the classes $C\{M_n\}$ and $C\{\frac{M_n}{M_0}\}$ coincide.\\

Using the convex regularization by means of the logarithm and the observations above, it is possible to show the Denjoy-Carleman theorem in its most general form (see \cite{Coh68} for a simple but detailed proof).
\begin{theorem}[The Denjoy-Carleman Theorem] \label{DJ-C}\ \\
Let $(M_n)_{n\in\NN_0}$ be a sequence of positive real numbers. Then the following conditions are equivalent
\begin{enumerate}
\item $C\{M_n\}$ is quasi-analytic,
\item $\sum\limits_{n=1}^\infty \frac{1}{\beta_n}=\infty$ with $\beta_n:=\inf_{k\geq n}\sqrt[k]{M_k},$
\item $\sum\limits_{n=1}^\infty \frac{1}{\sqrt[n]{M_n^c}}=\infty$,
\item $\sum\limits_{n=1}^\infty \frac{M_{n-1}^c}{M_n^c}=\infty$,
\end{enumerate}
where $({M}^c_n)_{n\in\NN_0}$ is the convex regularization of $({M}_n)_{n\in\NN_0}$ by means of the logarithm. 
\end{theorem}

Let us now state a simple result which has been repeatedly used throughout this paper.

\begin{prop}\label{ConstTimesCarl}\ \\
Let $({M_n})_{n\in\NN_0}$ be a sequence of positive real numbers. Then, $C\{M_n\}$ is quasi-analytic if and only if for any positive constant $\delta$ the class $C\{\delta M_n\}$ is quasi-analytic.\end{prop}
%\proof\ \\
%In the case $\delta=1$ the theorem trivially holds.\\
%
%\noindent Assume that $\sum\limits_{n=1}^{\infty}\frac{1}{\sqrt[n]{M_{n}}}=\infty$ and let us define $t(n):=\delta^{\frac 1{n}}$. 
%\begin{itemize} 
%\item If $0<\delta<1$ then $t(n)$ is increasing and so for any $n\in\NN$ we have $\delta^{\frac 1{n}} M_n^{\frac 1{n}}\leq\left(\lim\limits_{n\to\infty} \delta^{\frac 1{n}}\right)\, M_n^{\frac 1{n}}$. This implies that
%$$\sum_{n=1}^{\infty}\frac{1}{\sqrt[n]{\delta\, M_{n}}}\geq\sum_{n=1}^{\infty}\frac{1}{\sqrt[n]{M_{n}}}=\infty.$$
%\item If $\delta>1$ then $t(n)$ is decreasing and so for any $n\in\NN$ we have $\delta^{\frac 1{n}}\, M_n^{\frac 1{n}}\leq \delta\, M_n^{\frac 1{n}}$. This implies that
%$$\sum_{n=1}^{\infty}\frac{1}{\sqrt[n]{\delta\,M_{n}}}\geq\frac 1{{\delta}}\sum_{n=1}^{\infty}\frac{1}{\sqrt[n]{M_{n}}}=\infty.$$
% \end{itemize}
% 
%Assume that for any $\delta>0$ we have $\sum\limits_{n=1}^{\infty}\frac{1}{\sqrt[n]{\delta\, M_{n}}}=\infty$. Since $t(n)\to 1$ as $n\to\infty$, there exists $N\in\NN$ such that, for any $n\geq N$, we have that $t(n)\geq\frac 12$ and then
%$$\frac{1}{\sqrt[n]{\delta\, M_{n}}}\leq\frac{2}{\sqrt[n]{M_{n}}}.$$
%\endproof

In conclusion, let us introduce some interesting properties of log-convex sequences.
\begin{remark}\ \\
For a sequence of positive real numbers $(M_n)_{n\in\NN_0}$  the following properties are equivalent
\begin{description}
\item[(a)] $(M_n)_{n=0}^\infty$ is log-convex.
\item[(b)] $\left(\frac{M_{n}}{M_{n-1}}\right)_{n=1}^{\infty}$ is monotone increasing.
\item[(c)] $\left(\ln(M_n)\right)_{n=1}^{\infty}$ is convex.
\end{description}
\end{remark}
Note that the log-convexity is a necessary condition for a sequence to be a moment sequence.

\begin{prop}\label{corInc}\ \\
If the sequence $(M_n)_{n\in\NN_0}$ is log-convex and $M_0=1$, then $(\sqrt[n]{M_n})_{n=1}^{\infty}$ is monotone increasing.
\end{prop}
%\begin{proof}\ \\
%From (b) it follows that
%\[
%M_n = \frac{M_n}{M_0} = \prod_{j=1}^n \frac{M_j}{M_{j-1}} \leq
%\left(\frac{M_n}{M_{n-1}} \right)^n,
%\]
%i.e.
%\[
%M_{n-1}^{n} \leq M_{n}^{n-1},
%\]
%or equivalently
%\[
%M_{n-1}^{1/n-1} \leq M_{n}^{1/n}.
%\]
%\end{proof}

\begin{lemma}\label{CarlemanPari}\ \\
Assume that $(M_n)_{n\in\NN_0}$ is a log-convex sequence. The class $C\{M_n\}$ is quasi-analytic if and only if for any $j\in\NN$ the class $C\{\sqrt[j]{M_{jn}}\}$ is quasi-analytic.
%\begin{itemize}
%\item If $\sum\limits_{n=1}^\infty \frac{1}{\sqrt[n]{M_{n}}}=\infty$
%then for any $j\in\NN$
%$\sum\limits_{n=1}^\infty \frac{1}{\sqrt[jn]{M_{jn}}}=\infty.$
%\item If $\sum\limits_{n=1}^\infty \frac{1}{\sqrt[jn]{M_{jn}}}=\infty$ for some $j\in\NN$, then $\sum\limits_{n=1}^\infty \frac{1}{\sqrt[n]{M_{n}}}=\infty$.
%\end{itemize}
\end{lemma}
\proof\ \\
W.l.o.g.\!\! we can assume that $M_0=1$. (In fact, if $M_0\neq 1$ then one can always apply the following proof to the sequence $(\frac{M_n}{M_0})_{n\in\NN_0}$ by Proposition~\ref{ConstTimesCarl}.)
Let us first note that by Theorem~\ref{DJ-C} it is enough to prove that $\sum\limits_{n=1}^\infty \frac{1}{\sqrt[n]{M_{n}}}=\infty$ if and only if for all $j\in\NN$, $\sum\limits_{n=1}^\infty \frac{1}{\sqrt[jn]{M_{jn}}}=\infty$. Let us fix $j\in\NN$, then
\begin{eqnarray*}
\sum_{n=1}^\infty \frac{1}{\sqrt[n]{M_{n}}}\!\!\!\!&=&\!\!\!\!\sum_{n=1}^\infty \left(
  \frac{1}{\sqrt[jn]{M_{jn}}} + \frac{1}{\sqrt[jn+1]{M_{jn+1}}} +
  \ldots + \frac{1}{\sqrt[jn+(j-1)]{M_{jn+j-1}}}\right) + \sum_{n=1}^{j-1} 
  \frac{1}{\sqrt[n]{M_{n}}} \\
  &\leq& j \sum_{n=1}^\infty
  \frac{1}{\sqrt[jn]{M_{jn}}} + \sum_{n=1}^{j-1}
  \frac{1}{\sqrt[n]{M_{n}}},
\end{eqnarray*}
where the last inequality is due to Proposition \ref{corInc}. Hence, if $\sum\limits_{n=1}^\infty \frac{1}{\sqrt[n]{M_{n}}}$ diverges then $\sum\limits_{n=1}^\infty\frac{1}{\sqrt[jn]{M_{jn}}}$ diverges as well. On the other hand, if the series $\sum\limits_{n=1}^\infty \frac{1}{\sqrt[jn]{M_{jn}}}$ diverges for some $j\in\NN$, then also $\sum\limits_{n=1}^\infty \frac{1}{\sqrt[n]{M_{n}}}$ diverges since the latter contains more summands.\\
\endproof

\subsection{Complements about the space $\C_c^\infty(\RR^d)$}\label{Sec-App2}\ \\
Let us recall the definition of the inductive topology on $\C_c^\infty(\RR^d)$ (see~\cite[Vol. I, Section~V.4]{ReSi75} for a more detailed account on this topic).
\begin{definition}\label{DefInd}\ \\
Let $(\Lambda_n)_{n\in\NN}$ be an increasing family of relatively compact open subsets of $\RR^d$ such that $\RR^d=\bigcup\limits_{n\in\NN}\Lambda_n$. Let us consider the space $\C_c^\infty(\overline{\Lambda_n})$ of all infinitely differentiable functions on $\RR^d$ with compact support contained in $\overline{\Lambda_n}$ and let us endow $\C_c^\infty(\overline{\Lambda_n})$ with the Frech\'et topology generated by the directed family of seminorms given by 
\begin{equation}\label{norm2}
\left\|\varphi\right\|_{\leq a}:=\sum_{|\beta|\leq a}\max_{\r\in\overline{\Lambda_n}}\left|D^\beta\varphi(\r)\right|.
\end{equation}
Then as sets $$\C_c^\infty(\RR^d)=\bigcup_{n\in\NN}\C_c^\infty(\overline{\Lambda_n}).$$ We denote by $\mathscr{D}_{ind}(\RR^d)$ the space $\C_c^\infty(\RR^d)$ endowed with the inductive limit topology $\tau_{ind}$ induced by this construction. 
\end{definition}
It is easy to see that the previous definition is independent of the choice of the $\Lambda_n$'s. \\

In Subsection~\ref{Sec-GenFunct}, we gave a construction due to Y. M. Berezansky that allows to write $\C_c^\infty(\RR^d)$ as projective limit of a family of weighted Sobolev space (see~Definition~\ref{TeoBerez}). Berezansky actually proved that Definition~\ref{TeoBerez} is equivalent to the following standard one (see\! \cite[Chapter I, Section 3.10]{B86} for more details).

\begin{definition}\label{DProj}\ \\
Let $I$ be as in Definition~\ref{TeoBerez}, i.e. the set of all $k=(k_1, k_2(\r))$ such that $k_1\in\NN_0$, $k_2\in\C^\infty(\RR^d)$ with $k_2(\r)\geq 1$ for all $\r\in\RR^d$.  For each $k\in I$, let us introduce a norm on $\C_c^\infty(\RR^d)$ by setting 
$$\left\|\varphi\right\|_{\mathscr{D}_k(\RR^d)}:=\max_{\r\in\RR^d}\left(k_2(\r)\sum_{|\beta|\leq k_1}\left |(D^{\beta}\varphi)(\r)\right |\right).$$
Denote by $\mathscr{D}_k(\RR^d)$ the completion of $\C_c^\infty(\RR^d)$ w.r.t. the norm $\left\|\cdot\right\|_{\mathscr{D}_k(\RR^d)}$. Then as sets 
$$\C_c^\infty(\RR^d)=\bigcap_{k\in I} \mathscr{D}_k(\RR^d).$$
We denote by $\mathscr{D}_{proj}(\RR^d)$ the space $\C_c^\infty(\RR^d)$ endowed with the projective limit topology $\tau_{proj}$ induced by this construction.
\end{definition}

Furthermore, as already mentioned, Berezansky showed that $$\mathscr{D}_{proj}(\RR^d)=\projlim\limits_{(k_1, k_2(\r))\in I}W_2^{k_1}(\RR^d, k_2(\r) d\r)$$ is a nuclear space (where $I$ is as in Definition~\ref{TeoBerez}). The nuclearity of $\mathscr{D}_{proj}(\RR^d)$ follows from the fact that the index set $I$ always fulfills the following condition.

\begin{definition}[Condition (D)]\label{CondIndex}\ \\
We say that the set $K_0\subseteq I$ satisfies \emph{Condition (D)} if:\\
``For any pair $k=(k_1, k_2(\r))\in K_0$ there exists $k'=(k'_1, k'_2(\r))\in K_0$ such that 
\begin{itemize}
\item $k'_1\geq k_1+l$ (where $l$ is the smallest integer greater than $\frac d2$) 
\item $k'_2(\r)\geq \left(\max\limits_{|\beta|\leq l}| (D^\beta q)(\r)|\right)^2$, $\forall\,\r\in\RR^d$, for some function $q(\r)\in\C^l(\RR^d)$ chosen such that $$q^2(\r)\geq k_2(\r),\,\forall\,\r\in\RR^d\quad\text{and}\quad\int_{\RR^d}\frac{k_2(\r)}{q^2(\r)}d\r<\infty.$$
\end{itemize}
Note that the function $q(\r)$ depends on $k_2(\r)$ and $d$.''
\end{definition}

Condition (D) is sufficient for $\projlim\limits_{(k_1, k_2(\r))\in K_0}W_2^{k_1}(\RR^d, k_2(\r) d\r)$ to be nuclear.

Let us give some concrete examples of classes $K_0$ which satisfy Condition~(D) in the case $d=1$.

\begin{example}\ \\
Let $K_0:=\{ (k_1, k_2(r)) \,|\, k_1\in\NN_0, k_2(r)=C(1+r^{2n}), n\in\NN, 1\leq c\in\RR\}$. \\ 
Let us fix a pair $k=(k_1, k_2(r))\in K_0$, namely we fix $k=(k_1, C(1+r^{2n}))$ for some $k_1\in\NN_0$, some $n\in\NN$ and some real constant $C\geq 1$. For the same fixed $n$ and $C$, we define the function $q(r):=(2C(1+r^{2n+2}))^\frac{1}{2}\in\C^\infty(\RR)$.\\
Then we have that $q^2(r)=2C(1+r^{2n+2})\geq k_2(r)$ for all $r\in\RR$ and
 $$\int_{\RR}\frac{k_2(r)}{q^2(r)}d r=\int_{\RR}\frac{1+r^{2n}}{2(1+r^{2n+2})}d r<\infty.$$
Hence, using the special form of $q(r)$, we get that  $$\forall r\in\RR, \quad|Dq(r)|\leq (n+1)|q(r)|.$$  
Consequently, choosing $k'=(k'_1, k'_2(r))\in K_0$ such that 
$$k'_1:= k_1+1,\qquad k'_2(r):=(n+1)^2 q(r)^2,\quad\forall r\in\RR,$$  
we obtain that for all $r\in\RR$, $k'_2(r)\geq \left(\max\{|q(r)|,|Dq(r)|\}\right)^2$
and hence, Condition~(D) is fulfilled by $K_0$.
\end{example}

\begin{example}\ \\
Let $K_0:=\{ (k_1, k_2(r)) \,|\, k_1\in\NN_0, k_2(r)=1+e^{nr}, n\in\NN, 1\leq c\in\RR\}$. \\ 
Let us fix a pair $k=(k_1, k_2(r))\in K_0$, namely we fix $k=(k_1, C(1+e^{nr}))$ for some $k_1\in\NN_0$, some $n\in\NN$ and some real constant $C\geq 1$. For the same fixed $n$ and $C$, we define the function $q(r):=(C(1+e^{nr})(1+r^2))^\frac{1}{2}\in\C^\infty(\RR)$.\\Then we have that $q^2(r)=C(1+e^{nr})(1+r^2)\geq k_2(r)$ for all $r\in\RR$ and
 $$\int_{\RR}\frac{k_2(r)}{q^2(r)}d r=\int_{\RR}\frac{1}{1+r^2}d r<\infty.$$
Hence, using the special form of $q(r)$, we get that  $$\forall r\in\RR, \quad|Dq(r)|\leq \left(\frac{n}{2}+1\right)|q(r)|.$$ Consequently, if $B:=\sup\limits_{r\in\RR}\frac{(1+e^{nr})(1+r^2)}{1+e^{(n+1)r}}$ and if we choose $k'=(k'_1, k'_2(r))\in K_0$ such that 
 $$k'_1:= k_1+1,\qquad k'_2(r):=BC\left(\frac{n}{2}+1\right)^2(1+e^{(n+1)r}), \quad\forall r\in\RR,$$  
 then we obtain that for all $r\in\RR$, $$k'_2(r)\geq C\left(\frac{n}{2}+1\right)^2(1+e^{nr})(1+r^2)=\left(\frac{n}{2}+1\right)^2q^2(r)\geq\left(\max\{|q(r)|,|Dq(r)|\}\right)^2.$$
 \end{example}
 
\subsection{Construction of a total subset of test functions}\label{Sec-App3}\ \\
In this subsection, we provide an outline of the proof of Lemma~\ref{LemmaE} about the explicit construction of a set $E$ of the kind required in Definition~\ref{DefSeq}. For convenience, we give here the proofs only in the case when $E\subset\mathscr{D}_{proj}(\RR)$. The higher dimensional case follows straightforwardly.

For any $n\in\NN_0$, let $k^{(n)}:=(k_1^{(n)}, k_2^{(n)})\in I$, i.e.\! $k_1^{(n)}\in\mathbb{N}_0$ and $k_2^{(n)}: \mathbb{R} \rightarrow [1,\infty[$ such that $k_2^{(n)}\in\mathcal{C}^\infty(\RR)$. Let us consider the norm $\| \cdot\|_{H_{k^{(n)}}}$ defined in \eqref{NormW2weighted}, where $H_{k^{(n)}}:=W^{k_1^{(n)}}_2(\RR,k_2^{(n)}(x))$. We will denote by $\| \cdot\|_{H_{-k^{(n)}}}$ the norm on its dual space $W^{-k_1^{(n)}}_2(\RR,k_2^{(n)}(x))$.

Let $(d_n)_{n\in\NN_0}$ be a positive sequence which is not quasi-analytic, then there exists a non-negative infinite differentiable function $\varphi$ with support $[-1,1]$ such that for all $x\in\RR$ and $n\in\NN_0$ holds $|\frac{d^n}{dx^n}  \varphi(x)| \leq d_n$ (see \cite{Ru74}). Easy examples of increasing sequences of positive numbers which are not quasi-analytic are given by $n!(\ln n)^{2n}$ or $(n!)^{1+\varepsilon}$, for any $\varepsilon>0$.

\begin{lemma}\label{lemmino}\ \\
Let $(d_n)_{n\in\NN_0}$ be a log-convex increasing positive sequence which is not quasi-analytic, let $\varphi$ be as above. Define
\[
E_0 := \{ f_{y,p}(\cdot):=\varphi(\cdot -y) e^{ip\cdot}\ | \ y ,p\in \mathbb{Q}   \}.
\]
Then for any $y,p\in\mathbb{Q}$ and for any $n\in\NN_0$  we get
\[
 \| f_{y,p}\|_{H_{k^{(n)}}} \leq C_p^{k_1^{(n)}} d_{k_1^{(n)}} \sup_{x \in [-1,1]} \sqrt{ k_2^{(n)}(y+x)},
\]
where $C_p:= \sqrt{2}(1+|p|)$
and $E_0$ is total in $\mathscr{D}_{proj}(\mathbb{R})$.
\end{lemma} 

\proof\ \\
For any $y,p\in\mathbb{Q} $ we have that
\begin{align*}
  (\| f_{y,p}\|_{H_{k^{(n)}}})^2& \leq  \sum_{k=0}^{k_1^{(n)}} \int_{\mathbb{R}}\left( \sum_{l=0}^k \binom{k}{l} |p|^{k-l} \left| \frac{d^l}{dx^l}\varphi(x-y)\right| \right)^2  k_2^{(n)}(x) dx \\
%\\ & =  \sum_{k=0}^{k_1^{(n)}} \int_{\mathbb{R} } \sum_{l_1,l_2=0}^k \binom{k}{l_1}\binom{k}{l_2} |p|^{2k -l_1-l_2}
% \left| \frac{d^{l_1}}{dx^{l_1}}\varphi(x-y)\right|  \left| \frac{d^{l_2}}{dx^{l_2}}\varphi(x-y)\right|   k_2^{(n)}(x) dx 
 %\\ %& \leq \frac{1}{2} \sum_{k=0}^{k_1^{(n)}} \int_{\mathbb{R} } \sum_{l_1,l_2=0}^k \binom{k}{l_1}\binom{k}{l_2} |p|^{2k -l_1-l_2}2\left| \frac{d^{l_1}}{dx^{l_1}}\varphi(x-y)\right|^2    k_2^{(n)}(x) dx 
% \\ & \leq (1+|p|)^{k_1^{(n)}}\sum_{k=0}^{k_1^{(n)}} \int_{\mathbb{R} } \sum_{l_1=0}^k \binom{k}{l_1} |p|^{k -l_1}
% \left| \frac{d^{l_1}}{dx^{l_1}}\varphi(x-y)\right|^2     k_2^{(n)}(x) dx  \\ 
%& =  (1+|p|)^{k_1^{(n)}}\sum_{k=0}^{k_1^{(n)}}  \sum_{l=0}^k \binom{k}{l} |p|^{k -l}
 %\int_{\mathbb{R} } \left| \frac{d^{l}}{dx^{l}}\varphi(x-y)\right|^2    k_2^{(n)}(x) dx \\
  & \leq  (1+|p|)^{k_1^{(n)}}\sum_{k=0}^{k_1^{(n)}}  \sum_{l=0}^k \binom{k}{l} |p|^{k -l}
 \int_{[-1,1]} \left| \frac{d^{l}}{dx^{l}}\varphi(x)\right|^2    k_2^{(n)}(x+y) dx 
 \end{align*}
Using the bound for derivative of $\varphi$ and the fact that the sequence $(d_l)_l$ is monotone increasing
 we get
%\begin{align*}
%   (\| f_{y,p}\|_{H_{k^{(n)}}})^2
%%& \leq  2(1+|p|)^{k_1^{(n)}}\sum_{k=0}^{k_1^{(n)}}  \sum_{l=0}^k \binom{k}{l} |p|^{k -l}
%% (d_l)^2   \sup_{x \in[-1,1]} k_2^{(n)}(x+y)\\
% % &\leq 2(1+|p|)^{k_1^{(n)}} (d_{k_1^{(n)}})^2  (k_1^{(n)}+1)(1+ |p|)^{k_1^{(n)} } \sup_{x \in[-1,1]} k_2^{(n)}(x+y)\\
% &\leq   2^{k_1^{(n)}+1}(d_{k_1^{(n)}})^2 (1+|p|)^{2k_1^{(n)}} \sup_{x \in[-1,1]} k_2^{(n)}(x+y).
% \end{align*}
%and hence,
\begin{equation}\label{ineq0}
 \| f_{y,p}\|_{H_{k^{(n)}}}\leq\sqrt{2}\ d_{k_1^{(n)}}(\sqrt{2}(1+|p|))^{k_1^{(n)}}  \sup_{x \in[-1,1]} \sqrt{k_2^{(n)}(x+y)}.\end{equation}
 
Let us show that $E_0$ is total in $\mathscr{D}_{proj}(\mathbb{R})$.\\ 
If $E_0$ was not total, then by Hahn-Banach there would exist $\eta \in \mathscr{D}'_{proj}(\RR)$ with $\eta \neq 0$ such that for all $f \in span(E_0) $, $\eta(f) =0$. For such an $\eta$ we get in particular that
$\forall\ y,p\in\mathbb{Q},\,\langle f_{y,p}, \eta \rangle =0.$ Since the function $(y,p)\mapsto f_{y,p}$ from $\mathbb{Q}\times\mathbb{Q}$ to $\mathscr{D}_{proj}(\RR)$ is sequentially continuous, then 
\begin{equation}\label{ineq1}
\forall y,p\in\RR,\,\,\  \langle f_{y,p} , \eta \rangle =0.
\end{equation}\\ Let $\rho_\epsilon(\cdot):= \varepsilon^{-1}\rho(\varepsilon^{-1} \cdot )$ where $\rho$ is a non-negative function with compact support, i.e. $\rho_{\varepsilon}$ is an approximating identity  then
\begin{equation}\label{eq2}
\lim_{\varepsilon\downarrow 0} \int_{[-1,1]}f_{y,p}(x)  \rho_\varepsilon*\eta(x) dx
=\langle f_{y,p} , \eta \rangle=0,
\end{equation}
where the last equality is due to \eqref{ineq1}. Since $\eta$ is in some space $H_{-k^{(n)}}$ and as \eqref{ineq0}, holds, we get that
 \begin{equation} \label{eqqa27}
\left|\langle f_{y,p} , \rho_\varepsilon*\eta \rangle \right|
\leq \| f_{y,p}  \|_{H_{k^{(n)}}} \| \rho_\varepsilon*\eta\|_{H_{-k^{(n)}}}\leq c(1+|p|)^{k_1^{(n)}} \| \rho_\varepsilon*\eta\|_{H_{-k^{(n)}}},
\end{equation}
where $c:=d_{k_1^{(n)}}(\sqrt{2})^{k_1^{(n)}+1}\!\!\!\sup\limits_{x \in[-1,1]}\! \!\sqrt{k_2^{(n)}(x+y)}$ and so it depends only on $k_1^{(n)},k_2^{(n)},y$.\\
Since $\rho_\varepsilon$ is an approximating identity we get that  
\[
\lim_{\varepsilon\downarrow 0} \| \rho_\varepsilon*\eta\|_{H_{-k^{(n)}}} = \| \eta\|_{H_{-k^{(n)}}}
\]
The latter together with \eqref{eqqa27} imply that the function $\langle f_{y,p} , \rho_\varepsilon*\eta \rangle$ is uniformly bounded in $p$ and $\varepsilon$. By Lebesgue's dominated convergence theorem and by \eqref{eq2}, for any integrable function $\psi$ such that the Fourier transform $\hat{\psi} \in \mathscr{D}_{proj}(\mathbb{R})$ and for any $y\in\RR$ the following holds
%\begin{align*}
%0=& \lim_{\varepsilon\to 0} \int_{\mathbb{R}} \psi(p)\int_{[-1,1]} f_{y,p}(x)  \rho_\varepsilon*\eta(x) dx dp\\ 
%& \lim_{\varepsilon\to 0} \int_{\mathbb{R}} \psi(p)\int_{[-1,1]} \varphi(x-y) e^{ip x}  \rho_\varepsilon*\eta(x) dxdp \\ 
%& = \lim_{\varepsilon\to 0}\int_{[-1,1]} \varphi(x-y) \hat{\psi}(x)  \rho_\varepsilon*\eta(x) dx \\
%& = \lim_{\varepsilon\to 0} \langle \varphi(\cdot-y) \hat{\psi} , \rho_\varepsilon*\eta \rangle\\
%&=\langle \varphi(\cdot-y) \hat{\psi} , \eta \rangle = \langle  \hat{\psi} , \varphi(\cdot-y) \eta \rangle.
%\end{align*}
$$0= \lim_{\varepsilon\to 0} \int_{\mathbb{R}} \psi(p)\int_{[-1,1]} f_{y,p}(x)  \rho_\varepsilon*\eta(x) dx dp\langle \varphi(\cdot-y) \hat{\psi} , \eta \rangle = \langle  \hat{\psi}, \varphi(\cdot-y) \eta \rangle. $$
As any test-function in $\mathscr{D}_{proj}(\mathbb{R})$ is of the form  $\hat{\psi}$, %Fourier inversion theorem for S
we have that also as a distribution for any $y\in\RR$,
$
\varphi(\cdot-y) \eta\equiv0.
$ \\
Since $\varphi$ is not zero there exists an open ball $B$ on which $\varphi$ is never zero. Define a partition of unity $(\chi_n)_{n\in\NN_0}$, where each $\chi_n$ is supported in a ball of the form $y_n + B$.
Hence, for all $\psi \in\C_c^\infty(\mathbb{R})$ 
%\[
%\psi (\cdot)= \sum_{n=0}^\infty \chi_n(\cdot) \varphi(\cdot-{y_n}) \frac{\psi(\cdot)}{\varphi(\cdot-{y_n})}
%\]
%and so
\[
\langle \psi,  \eta \rangle = \sum_{n=0}^\infty \langle \chi_n(\cdot) \frac{\psi(\cdot)}{\varphi(\cdot-{y_n})}, \varphi(\cdot-{y_n}) \eta \rangle=0,
\]
which means that $\eta \equiv 0$.
\endproof
Making use of the previous result, we are going to prove Lemma~\ref{LemmaE} that we rewrite here for convenience.
\begin{lemma}\ \\
Let $(c_n)_{n\in\NN_0}$ be an increasing sequence of positive numbers which is not quasi-analytic.
Then the set
\[
E := \left\{f \in \mathscr{D}_{proj}(\mathbb{R}) \left| \forall\ n\in\NN_0,\  \| f\|_{H_{k^{(n)}}} \leq c_{k_1^{(n)}}  \sup_{\stackrel{z\in\RR}{|z|\leq n}}\sup_{x \in [-1,1]}\sqrt{ k_2^{(n)}(z+x)} \right.\right\}
\]
contains a countable subset which is total in $\mathscr{D}_{proj}(\mathbb{R})$. Hence, $E$ is total in $\mathscr{D}_{proj}(\mathbb{R})$.
\end{lemma}

\proof\ \\
Let us first show that the proof reduces to find an increasing sequence $(d_n)_{n\in\NN_0}$
of positive numbers which is not quasi-analytic and which is such that for any real constant $C >0$ 
\begin{equation}\label{rel-asintotica}
\lim_{j \rightarrow \infty} \frac{C^j d_j}{c_j}=0.
\end{equation}  
%In fact, if $\lim\limits_{n\to\infty}k_1^{(n)}=\infty$ then \eqref{rel-asintotica} guarantees that $\sup\limits_{n} \frac{C^{k_1^{(n)}}d_{k_1^{(n)}}}{c_{k_1^{(n)}}}<\infty$ and if instead $\lim\limits_{n\to\infty}k_1^{(n)}=l<\infty$ then $\sup\limits_{n} \frac{C^{k_1^{(n)}}d_{k_1^{(n)}}}{c_{k_1^{(n)}}}=\frac{C^{l}d_{l}}{c_{l}}<\infty$. Hence, in both cases we can always define 
In this case, we can always define
$
\frac{1}{q} := \sup\limits_{n}\frac{C^{k_1^{(n)}}d_{k_1^{(n)}}}{c_{k_1^{(n)}}},
$
and so, by Lemma~\ref{lemmino}, for any $y,p\in\mathbb{Q}$, every function of the form $qf_{y,p}$ is such that
\[
 \| q f_{y,p}\|_{H_{k^{(n)}}}\leq q C_p^{k_1^{(n)}} d_{k_1^{(n)}} \sup_{x \in [-1,1]} \sqrt{ k_2^{(n)}(y+x)}\leq c_{k_1^{(n)}} \sup_{\stackrel{z\in\RR}{|z|\leq n}}\sup_{x \in [-1,1]} \sqrt{ k_2^{(n)}(y+x)}.
\]
Hence, the set $E$ contains $qE_0$. Consequently, since $E_0$ is  total in $\mathcal{D}(\RR^d)$, the same is true for $qE_0$ and hence, for $E$. \\
It remains to construct an increasing sequence $(d_n)_n$ of positive numbers not quasi-analytic and such that \eqref{rel-asintotica} holds. First note that our requirement is equivalent to define an increasing sequence $(d_n)_n$ of positive numbers such that  $\sum_{n=1}^\infty\frac{1}{\sqrt[n]{d_n}}<\infty$ and
$
\lim_{n \rightarrow \infty } \frac{\sqrt[n]{d_n}}{\sqrt[n]{c_n}}=0.
$
Indeed, for each $C$ and for each $\varepsilon>0$ there exists $N$ such that for all $n \geq N$ holds $d_n \leq \left(\frac{\varepsilon}{C}\right)^n  c_n$ and hence also $C^n d_n \leq \varepsilon^n c_n$. \\ 
Our problem reduces to find, given a decreasing sequences $(a_n)_n$ of positive numbers with $\sum_{n=1}^\infty a_n<\infty$, a decreasing sequence $ (b_n)_n$ of positive numbers such that $\sum_{n=1}^\infty b_n<\infty$ and
$
\lim_{n \rightarrow \infty} \frac{b_n}{a_n} =\infty.
$\\
For any $k\in\NN$ let us define $N_k := \min\{ m | \sum_{n=m}^\infty a_n \leq \frac{1}{k^2} \}$ and also
\[
b_n := \min\left\{a_n \left( 1 + \sum_{k\in \mathbb{N} \ : \ N_k \leq n } \sqrt{k} \right), b_{n-1}\right\}, 
\]
with $b_0:=a_0 \left( 1 + \sum\limits_{k\in \mathbb{N} \ : \ N_k = 0 } \sqrt{k} \right)$.
Then
\begin{eqnarray*}
\sum_{n=1}^\infty b_n &\leq&\sum_{n=1}^\infty a_n\left( 1 + \sum_{k\in \mathbb{N} \ : \ N_k \leq n } \sqrt{k} \right)
%&\leq&\sum_{n=1}^\infty a_n+\sum_{n=1}^\infty a_n \sum_{k\in \mathbb{N} \ : \ N_k \leq n } \sqrt{k} \\
%&\leq&\sum_{n=1}^\infty a_n+ \sum_{k=1}^\infty  \sqrt{k} \sum_{n =N_k}^\infty a_n \\
\leq\sum_{n=1}^\infty a_n+ \sum_{k=1}^\infty k^{-3/2} < \infty,
\end{eqnarray*}
It follows that $\lim_{n\to\infty} b_n=0$. Then latter together with the definition $(b_n)_n$ implies that there exists an infinite subsequence $(b_{n_j})_{j}\subset(b_n)_n$ such that
$$
\forall\, j\in\NN \,:\, b_{n_j}=a_{n_j} \left( 1 + \sum_{k\in \mathbb{N} \ : \ N_k \leq n_j } \sqrt{k} \right).
$$
%In fact, if there was only a finite number of $n$ such that $b_n=a_{n} \left( 1 + \sum_{k\in \mathbb{N} \ : \ N_k \leq n} \sqrt{k} \right)$ then $\lim_{n\to\infty} b_n=0$ would have implied that $\exists\, n_1\in\NN\:\ \forall\ n\geq n_1\ b_n\equiv 0$,  which is impossible since for all $n\in\NN_0$, $b_n>0$.
For such a subsequence we have that
%Since $(N_k)_k$ is increasing and $\lim\limits_{k \rightarrow \infty}N_k=\infty$, 
\begin{equation}\label{subseq}
\lim_{j\to\infty}\frac{b_{n_j}}{a_{n_j}}=\lim_{j\to\infty}\left( 1 + \sum_{k\in \mathbb{N} \ : \ N_k \leq n_j } \sqrt{k} \right)=\left( 1 + \sum_{k=1}^\infty \sqrt{k} \right)=\infty.
\end{equation}
Now let us note that for any $n\in\NN$ we have either that
$\frac{b_n}{a_n}=\frac{b_{n-1}}{a_n}\geq \frac{b_{n-1}}{a_{n-1}}$
or that
$$\frac{b_n}{a_n}=\frac{a_{n} \left(1 +\!\!\!\!\!\!\!\!\!\sum\limits_{k\in \mathbb{N} \ : \ N_k \leq n }\!\!\!\sqrt{k} \right)}{a_n}=\left(1 +\!\!\!\!\!\!\!\!\! \sum\limits_{k\in \mathbb{N} \ : \ N_k \leq n }\!\!\!\sqrt{k} \right)\geq\left(1 +\!\!\!\!\!\!\!\!\!\sum\limits_{k\in \mathbb{N} \ : \ N_k \leq {n-1} }\!\!\! \sqrt{k} \right)\geq \frac{b_{n-1}}{a_{n-1}}.$$
Hence, the sequence $({b_n}/{a_n})_n$ is increasing and has a subsequence such that \eqref{subseq} holds, then we get that
$\lim_{n\to\infty}\frac{b_n}{a_n}=\infty.$\\
\endproof

\section*{Acknowledgments}
The authors would like to thank Yuri Kondratiev, Eugene Lytvynov and Martin Kolb for helpful discussions.

\bibliographystyle{plain}

\begin{thebibliography}{10}

\bibitem{B86}
Y.~M. Berezansky.
\newblock {\em Selfadjoint Operators in Spaces of Functions of Infinite Many
  Variables}, volume~63 of {\em Trans.~Amer.~Math.~Soc.}
\newblock American Mathematical Society, 1986.

\bibitem{BeKo88}
Y.~M. Berezansky and Y.~G. Kondratiev.
\newblock {\em Spectral Methods in Infinite-Dimensional Analysis}.
\newblock Naukova Dumka, Kiev, 1988.
\newblock In Russian. English translation: {K}luwer {A}cademic {P}ublishers,
  {D}ordrecht, 1995.

\bibitem{BeKoKuLy99}
Y.~M. Berezansky, Y.~G. Kondratiev, T.~Kuna, and E.~W. Lytvynov.
\newblock On spectral representation for correlation measures in configuration
  space analysis.
\newblock {\em Methods. Funct. Anal. and Topology}, 5:87--100, 1999.

\bibitem{BS71}
Y.~M. Berezansky and S.~N. {\v{S}}ifrin.
\newblock The generalized degree symmetric moment problem.
\newblock {\em Ukrain. Mat. \v Z.}, 23:291--306, 1971.

\bibitem{Berg87}
C.~Berg.
\newblock The multidimensional moment problem and semi-groups.
\newblock {\em {P}roc. {S}ymp. {A}ppl. {M}ath.}, 37:110--124, 1987.

\bibitem{BeMa}
C.~Berg and P.~H. Maserick.
\newblock Polynomially positive definite sequences.
\newblock {\em Mathematische Annalen}, 259:487--495, 1982.

\bibitem{BKM04}
E.~N. Brown, R.~E. Kass, and P.~P. Mitra.
\newblock Multiple neural spike train data analysis: state-of-the-art and
  future challenges.
\newblock {\em Nature Neuroscience}, 7:456--471, 2004.

\bibitem{Carl26}
T.~Carleman.
\newblock {\em Les fonctions quasi-analytiques}, volume~7 of {\em Collection de
  monographies sur la th{\'e}orie des fonctions publi{\'e}e sous la direction
  de M. E. Borel.}
\newblock Gauthier-Villars, Paris, 1926.

\bibitem{CimMarNet11}
J.~Cimpri{\v{c}}, M.~Marshall, and T.~Netzer.
\newblock Closures of quadratic modules.
\newblock {\em Israel J. Math.}, 183, 2011.

\bibitem{Coh68}
P.~J. Cohen.
\newblock A simple proof of the {D}enjoy-{C}arleman theorem.
\newblock {\em The American Mathematical Monthly}, 75(1):26--31, 1968.

\bibitem{DeJeu03}
M. de~Jeu.
\newblock Determinate multidimensional measures, the extended {C}arleman
  theorem and quasi-analytic weights.
\newblock {\em Ann. Probab.}, 31(3):1205--1227, 2003.

\bibitem{Devi53}
A.~Devinatz.
\newblock Integral representations of positive definite functions.
\newblock {\em Trans. Amer. Math. Soc.}, 74:56--77, 1953.

\bibitem{Gel-Shi68}
I.~M. Gel{\cprime}fand and G.~E. Shilov.
\newblock {\em Generalized functions. {V}ol. 2}.
\newblock Academic Press [Harcourt Brace Jovanovich Publishers], New York, 1968
  [1977].
\newblock Spaces of fundamental and generalized functions, Translated from Russian by Morris D. Friedman, Amiel Feinstein and Christian P. Peltzer.

\bibitem{HaMcDo87}
J.~P. Hansen and I.~R. McDonald.
\newblock {\em Theory of simple liquids}.
\newblock Academic Press, New York, 2nd edition, 1987.

\bibitem{Hav36}
E.~K. Haviland.
\newblock On the moment problem for distribution functions in more than one
  dimension {II}.
\newblock {\em Amer. J. Math.}, 58:164--168, 1936.

\bibitem{JM04}
M.~R. Jarvis and P.~P. Mitra.
\newblock Sampling properties of the spectrum and coherency of sequences of
  action potentials.
\newblock {\em Neural Comp.}, 13:717--749, 2004.

\bibitem{KoKuOl02}
{Y}.~G. Kondratiev, T.~Kuna, and M.~J. Oliveira.
\newblock Holomorphic {B}ogoliubov functionals for interacting particle systems
  in continuum.
\newblock {\em J. Funct. Anal.}, 238(2):375--404, 2006.

\bibitem{KuMarSch05}
S.~Kuhlmann, M.~Marshall, and N.~Schwartz.
\newblock Positivity, sums of squares and the multi-dimensional moment problem.
  {II}.
\newblock {\em Adv. Geom.}, 5(4):583--606, 2005.

\bibitem{Ku09}
T.~{Kuna}, J.~L. {Lebowitz}, and E.~R. {Speer}.
\newblock {Necessary and sufficient conditions for realizability of point
  processes}.
\newblock {\em The Annals of Applied Probability}, 21(4):1253--1281, 2011.

\bibitem{MolLach}
R.~Lachieze-Rey and I.~Molchanov.
\newblock Regularity conditions in the realisability problem in applications to
  point processes and random closed sets.
\newblock {\em To appear in Annals of Applied Probability.}

\bibitem{LasBook}
J.~B. Lasserre.
\newblock {\em Moments, positive polynomials and their applications}, volume~1
  of {\em Imperial College Press Optimization Series}.
\newblock Imperial College Press, London, 2010.

\bibitem{Las2011}
J.~B. Lasserre.
\newblock The {$\bold{K}$}-moment problem for continuous linear functionals.
\newblock {\em Trans. Amer. Math. Soc.}, 365(5):2489--2504, 2013.

\bibitem{LauMo}
M.~Laurent.
\newblock Sums of squares, moment matrices and optimization over polynomials.
\newblock {\em Emerging Application of Algebraic Geometry}, 149:157--270, 2008.

\bibitem{Le75a}
A.~Lenard.
\newblock States of classical statistical mechanical systems of infinitely many
  particles {I}.
\newblock {\em Arch. Rational Mech. Anal.}, 59:219--239, 1975.

\bibitem{Le75b}
A.~Lenard.
\newblock States of classical statistical mechanical systems of infinitely many
  particles {II}.
\newblock {\em Arch. Rational Mech. Anal.}, 59:241--256, 1975.

\bibitem{Mand52}
S.~Mandelbrojt.
\newblock {\em S\'eries adh\'erentes, r\'egularisation des suites,
  applications}.
\newblock Gauthier-Villars, Paris, 1952.

\bibitem{MarshBook}
M.~Marshall.
\newblock {\em Positive polynomials and sums of squares}, volume 146 of {\em
  Mathematical Surveys and Monographs}.
\newblock American Mathematical Society, Providence, RI, 2008.

\bibitem{Me77}
J.~Mecke.
\newblock Eine {C}arakterisierung des {W}estcottschen {F}unktionals.
\newblock {\em Math. Nachr.}, 80:295--313, 1977.

\bibitem{Mol05}
I.~S. Molchanov.
\newblock {\em Theory of random sets}.
\newblock Probability and its applications. Springer, New York, 2005.

\bibitem{MDL04}
D.~J. Murrell, U.~Dieckmann, and R.~Law.
\newblock On moment closure for population dynamics in continuous space.
\newblock {\em J. Theo. Bio.}, 229:421--432, 2004.

\bibitem{Nuss65}
A.~E. Nussbaum.
\newblock Quasi-analytic vectors.
\newblock {\em Ark. Mat.}, 6:179--191, 1965.

\bibitem{PolSze76}
G.~P{\'o}lya and G.~Szego.
\newblock {\em Problems and theorems in analysis. Vol. II}.
\newblock Springer Verlag, 1976.

\bibitem{PowSch01}
V.~Powers and C.~Scheiderer.
\newblock The moment problem for non-compact semialgebraic sets.
\newblock {\em Adv. Geom.}, 1(1):71--88, 2001.

\bibitem{Put93}
M.~Putinar.
\newblock Positive polynomials on compact semi-algebraic sets.
\newblock {\em Indiana Univ. Math. J.}, 42(3):969--984, 1993.

\bibitem{ReSi75}
M.~Reed and B.~Simon.
\newblock {\em Methods of modern mathematical physics}, vols. I and II.
\newblock Academic Press, New York and London, 1975.

\bibitem{Riesz23}
M.~Riesz.
\newblock {\em Sur le probl{\`e}me de moments: Troisi{\`e}me note}, volume~17.
\newblock Arkiv f{\"o}r matematik, Astronomi och Fysik, 1923.

\bibitem{Ru74}
W.~Rudin.
\newblock {\em Real and Complex Analysis}.
\newblock McGraw-Hill, New York, second edition, 1974.

\bibitem{Schm91}
K.~Schm\"{u}dgen.
\newblock The ${K}-$moment problem for compact semi-algebraic sets.
\newblock {\em Math. Ann.}, 289:203--206, 1991.

\bibitem{Sch57}
L.~Schwartz.
\newblock {\em Th\'eorie des distributions}, volume~1.
\newblock Hermann, 1957.

\bibitem{Sch73}
L.~Schwartz.
\newblock {\em Radon measures on arbitrary topological spaces and cylindrical
  measures}.
\newblock Published for the Tata Institute of Fundamental Research, Bombay by
  Oxford University Press, London, 1973.

\bibitem{Sh-Tam43}
J.~A. Shohat and J.~D. Tamarkin.
\newblock {\em The {P}roblem of {M}oments}.
\newblock American Mathematical Society Mathematical surveys, vol. I. American
  Mathematical Society, New York, 1943.

\bibitem{S74}
S.~N. {\v{S}}ifrin.
\newblock Infinite-dimensional symmetric analogues of the {S}tieltjes problem
  of moments.
\newblock {\em Ukrain. Mat. \v Z.}, 26:696--701, 718, 1974.

\bibitem{Stoy00}
D.~Stoyan.
\newblock Basic ideas of spatial statistics.
\newblock In {\em Statistical physics and spatial statistics ({W}uppertal,
  1999)}, volume 554 of {\em Lecture Notes in Phys.}, pages 3--21. Springer,
  Berlin, 2000.

\bibitem{To02}
S.~Torquato.
\newblock {\em Random Heterogeneous Materials: Macroscopic Properties}.
\newblock Springer-Verlag, Berlin, Heidelberg, NewYork, 2002.

\bibitem{Tre67}
F.~Tr{\`e}ves.
\newblock {\em Topological vector spaces, distributions and kernels}.
\newblock Academic Press, New York, 1967.

\bibitem{Ze83}
H.~Zessin.
\newblock The method of moments for random measures.
\newblock {\em Z. Wahrsch. Verw. Gebiete}, 62(3):395--409, 1983.

\end{thebibliography}
\def\cprime{$'$}

\end{document}